\documentclass[preprint,pra,aps,showpacs,showkeys]{revtex4}
\usepackage{graphicx}
\usepackage{amssymb}
 
\begin{document}
 
\title{Numerical Representation of the Incomplete Gamma Function of Complex Argument}
 
\author{RICHARD J. MATHAR}
\email{mathar@mpia.de}
\homepage{http://www.mpia.de/MIDI/People/mathar}
\affiliation{
Max-Planck Institute of Astronomy, K\"onigstuhl 17, 69117 Heidelberg, Germany}

\def\erf{\mathop{\rm erf}\nolimits}

\date{\today}
 
\begin{abstract}
Various approaches to the numerical representation of the Incomplete
Gamma Function $F_m(z)$ for complex arguments $z$ and small integer indexes $m$
are compared with respect to numerical fitness (accuracy and speed). We consider
power series, Laurent series,
Gautschi's 
approximation to the Faddeeva function, classical numerical methods of
treating the standard integral representation,
and others not yet covered by the literature.

The most suitable scheme is the construction of Taylor expansions around
nodes of a regular, fixed grid in the $z$-plane, which stores
a static matrix of higher derivatives. This is the
obvious extension to a procedure often in use for real-valued $z$.

\end{abstract}
\maketitle

\section{Overview}
\subsection{Motivation}
The Incomplete Gamma Function $F_m(z)$ is at the heart of the computation of 
Electron Repulsion Integrals over Gaussian-type basis functions
\cite{Boys,Zivko,Schwerdtfeger1988,Obara}.
If these are attached to moving atoms, and the dominant part of
the time dependency is kept with the bases (instead of being hidden in the
expansion coefficients), the argument $z$ is complex-valued (App.\ \ref{app.travo}).

With respect to the index $m$, we look at this as being derived
from coupling of integer-valued orbital quantum numbers, whence
deal only with small, non-negative $m$ unless otherwise noted.

Several generalizations \cite{Kilbas,Paris2003,Miller} are not considered here,
except the---important---identification with some confluent
hypergeometric series. 

\subsection{Contents}
The explicit intention of this script is to compare a wider range of
methods than proposed on the same subject before \cite{Carsky}.
Continued fractions are not covered here, because they have already been
detailed before \cite{Jones,Luke1975Book}.
Also, the Temme and Paris approximations \cite{Temme,Paris2002} have been left
aside, because
they would work around the Complementary Error Function of a complex argument,
the calculation of which is already of the same order of complexity as
the original problem.
No attempt is made to discuss the Gauss-Rys quadrature \cite{King,Dupuis2001,Lindh,Ishida1991}
for complex $z$, which constructs a system of orthogonal polynomials over
a finite interval with weight function $\exp(-zt^2)$;
the explicit notation of its polynomial of degree 1, Eq.\ (2.4)
by Steen et al.\ \cite{Steen}, illustrates that this also would
start from the complex Error Function. 
Barakat \cite{Barakat}\cite[(13.3.9)]{AS} reports on $z$-values on the imaginary axis mapped on
Bessel functions with real-valued arguments, which we do not follow
on the same reasoning.

For the case of real-valued $z$ and large $m$ we refer the reader to the
article by Takenaga \cite{Takenaga}.

$F_m(z)$ is eventually defined by an integral with a rather simple
kernel. The following chapters are roughly grouped according to how
much effort is spent on isolating special
aspects of this kernel in exactly integrable terms, in hope of 
catching up the oscillations induced by $\Im z$ or the steep slopes
induced for large $\Re z$, for example.

\subsection{Fundamentals}
We define $F_m(z)$ through its simplest integral representation,
\begin{equation}
F_m(z)\equiv\int_0^1 t^{2m}e^{-zt^2}dt=\frac{1}{2}\int_0^1 u^{m-1/2}e^{-zu}du .
\label{eq.FmInt}
\end{equation}

Though the application in the literature often focuses
on the positive real axis (as $z=x+iy$ represents some square of real-valued
widths of orbital exponentials then),
we will cover the general case.

Anyway, the complex-conjugate symmetry
\begin{equation}
F_m(\overline z)=\overline{F_m(z)},
\end{equation}
(with $\overline z\equiv \Re z-i\Im z$) allows us to
restrict the analysis to the
cases of $\Im z\ge 0$.

The forward recurrence \cite{Takashima}\cite[(6.5.21)]{AS}
\begin{equation}
2zF_m(z)=(2m-1)F_{m-1}(z)-e^{-z} \qquad (m>\frac{1}{2})
\label{eq.recurr}
\end{equation}
is a rephrasing of \cite[(13.4.4)]{AS}.
It suffers from cancellation of about $-\log_{10}|2z|$ ($m=0$) and 
$-\log_{10}|z/m|$ ($m>0$) decimal places when applied to $|z|\ll 1$ \cite{Murchie}.
The corresponding backward recursion suffers from cancellation of  about
$-\log_{10}|(m-1/2)/z|$ decimal places,
if $\Re z$ is negative and $|z|$ large.
Generally speaking, the recursion allows to pick a numerically favorable $m$ and switch
to others at a small additional cost.

Further down, numerical precision is demonstrated by the number of decimal
places, defined as the negative Brigg Logarithm of the relative error,
$d\equiv-\log_{10}|1-\hat F_m(z)/F_m(z)|$, given a high precision
actual value $F_m(z)$ and its approximation $\hat F_m(z)$, both computed
with Maple at 30 decimal places.

\section{Hypergeometric Series}\label{sec.referenc}
Expansion of the exponential in (\ref{eq.FmInt})
and term-by-term integration yields the
confluent hypergeometric series (Kummer Function)
\begin{equation}
F_m(z)=\frac{1}{2m+1}{}_1F_1(m+\frac{1}{2};m+\frac{3}{2};-z); \qquad (m>-1/2)
\label{eq.Kumm}
\end{equation}
\begin{equation}
_1F_1(a;a+1;z)=1+\frac{az}{a+1}+\cdots+\frac{az^n}{(a+n)n!}+\cdots
\label{eq.KummTay}
\end{equation}

\begin{figure}
\includegraphics[bb=144 191 496 820]{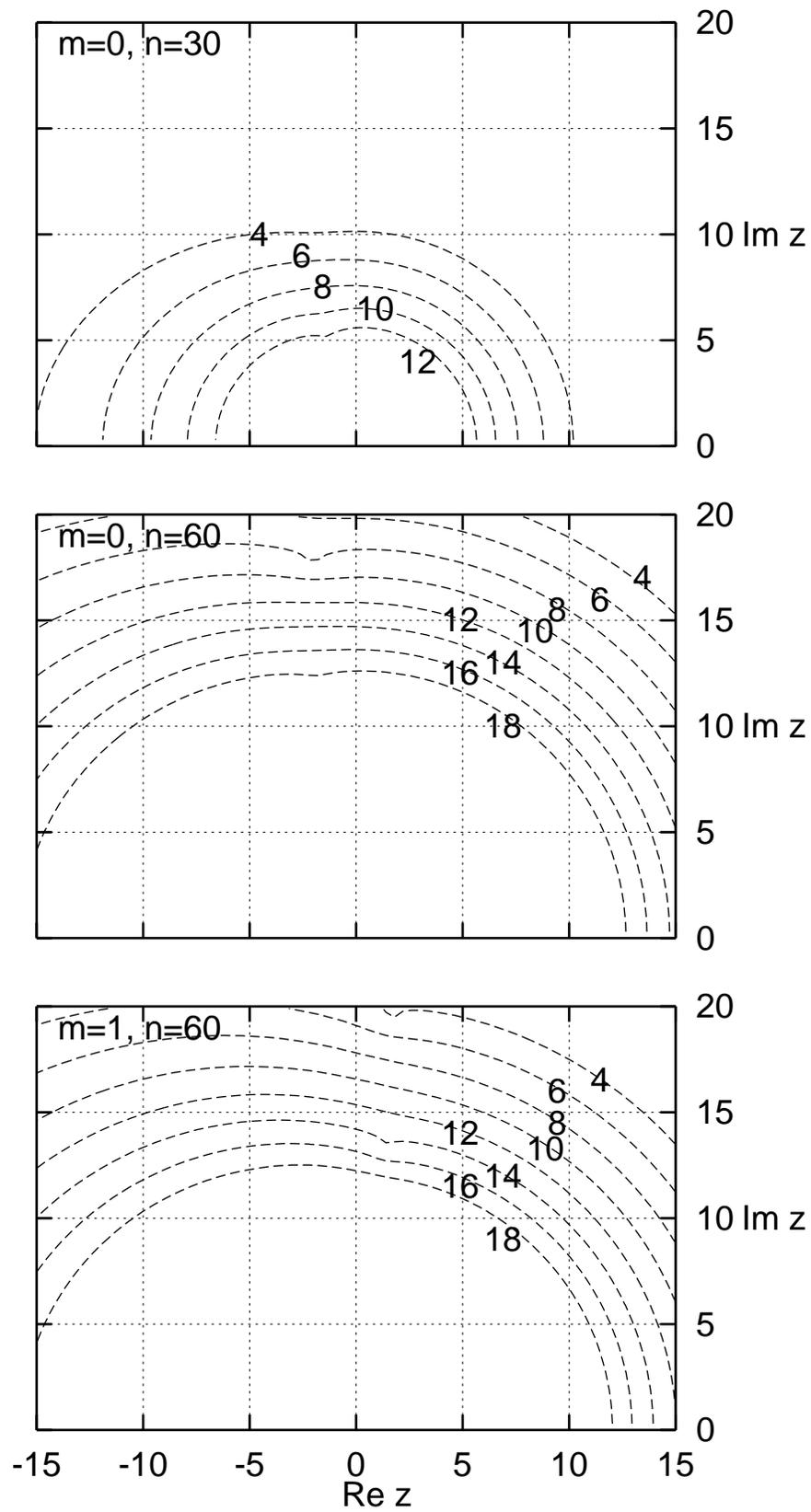}
\caption{
Contour levels of the number of valid decimal digits 
$d=4\ldots 18$ of $F_m(z)$ by the power series (\ref{eq.Kumm}) and (\ref{eq.KummTay}),
if the power series is truncated after $n=30$ or $60$ terms.
\label{fig.tayl}
}
\end{figure}

We start with the hypergeometric power series as a reference
because there is an evident and flexible implementation
with a clear numerical cost of one complex multiplication
(which could be implemented with
three real multiplications and five real-valued additions,
or with four real-valued multiplications and two real-valued additions
\cite[3.7.2]{Nussbaumer})
and addition per term:
one accumulates terms until the new term's contribution falls below a limit
set by a preset relative error.
Fig.\ \ref{fig.tayl} verifies
that the series converges best close to the origin of the complex plane,
as expected for any power series.
The unexpected feature is that the series performs worse if terms are 
alternating---that is if $z$ is close to the positive real axis in
Fig.\ \ref{fig.tayl}---than in the case of the non-alternating mirror
point $-\overline z$.
This is a by-product of a growth of the terms modulus up to the partial sum
of index $n\approx |z|-1$. The alternating case must overcome a massive
cancellation of digits when passing this index. Consequently it needs much
more terms until the partial sums approach the order of magnitude of the exact
result.
If measured in terms of the {\em absolute\/}
error after summation of $n$ terms, the alternating case would indeed
perform better.

The Shanks transformation \cite{Wynn1956} $\epsilon_2(S_n)$
of the partial sums $S_n$ of (\ref{eq.KummTay})---which is to this second order
just the Aitken transformation \cite[(3.9.7)]{AS}---
would improve the accuracy of the plots of Fig.\ \ref{fig.tayl}
by roughly 1.5 digits.
This involves handling of finite differences between numbers that are
(supposedly) already close to each other and is more tricky than the 
analytical transformation formulas indicate.

\section{Laurent Series}
\subsection{Barnes' Analysis}\label{sec.laur}

\begin{figure}
\includegraphics[bb=156 398 561 824]{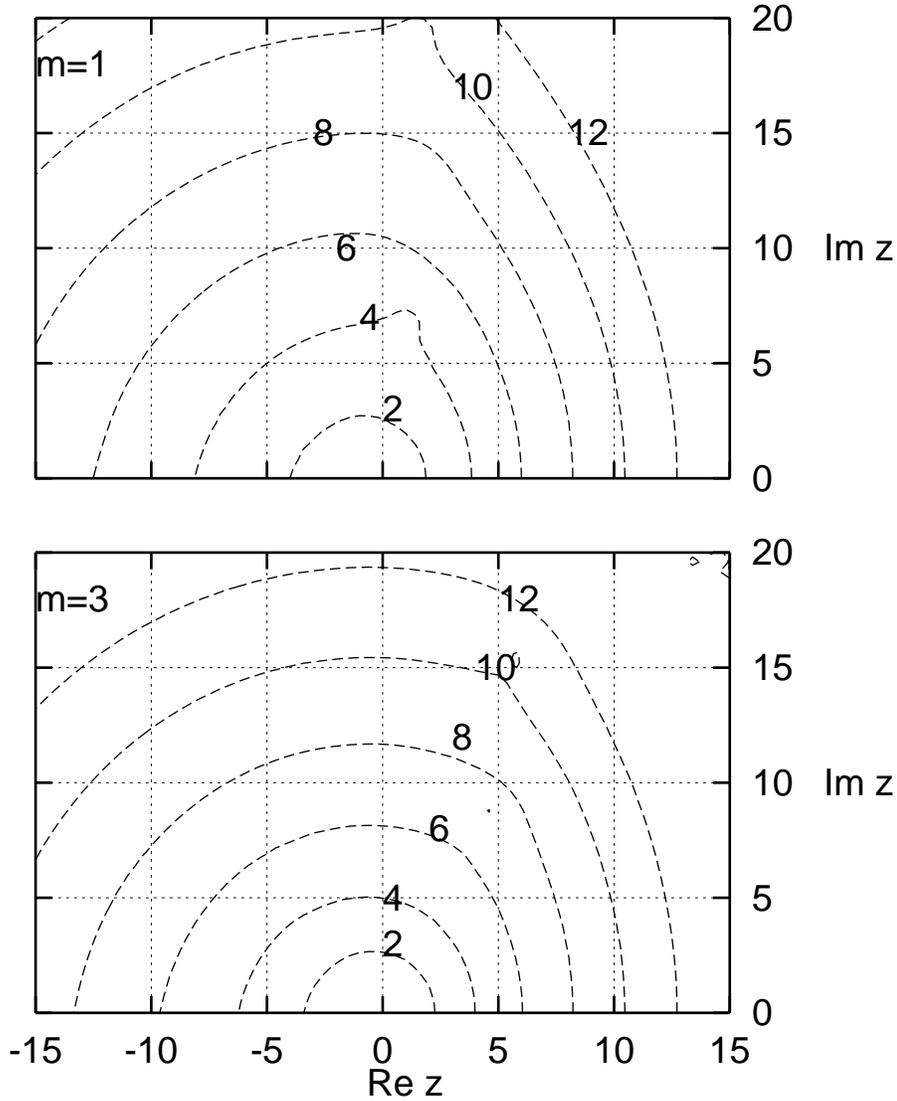}
\caption{
Contour levels of the number of valid decimal digits 
$d=2\ldots 12$ of $F_m(z)$, if the Laurent series (\ref{eq.KummLau})
is summed up to $n=|z|+a$.
\label{fig.lauren}
}
\end{figure}

The asymptotically convergent Laurent series for large $|z|$ is a special case
of Eq.\ (13.5.1) in \cite{AS} or taken from \S 6 in \cite{Barnes}:
\begin{equation}
_1F_1(a;a+1;z)=\Gamma(a+1)(-z)^{-a}+\frac{ae^z}{z}\sum_{n=0}^\infty (1-a)_n z^{-n},
\label{eq.KummLau}
\end{equation}
with Pochhammer's Symbol defined as \cite[(6.1.22)]{AS}
\begin{equation}
(b)_n\equiv \Gamma(b+n)/\Gamma(b)=\left\{\begin{array}{cc}
1 & , (n=0) \\
b(b+1)(b+2)\cdots (b+n-1) & , (n>0) \\
\end{array}\right.
\end{equation}
This series is also known under the label ``high-T'' expansion in quantum
chemistry \cite{Murchie,Obara}.
The Gamma Function is not of concern since it is only needed for
half-integer values, and would be tabulated based on \cite[(6.1.12)]{AS}.
Asymptotic convergence means that the terms in (\ref{eq.KummLau}) shrink until
$n\le |z|+a$ and grow afterwards. This inherent limitation to the
achievable
accuracy is put into concrete with Fig.\ \ref{fig.lauren}.

Swapping the the sign of $z$ in Eq.\ (\ref{eq.Kumm}),
the series becomes alternating near the positive real axes
of the plots in Fig.\ \ref{fig.lauren} for $\Re z>0$, which leads
to some obvious left-right asymmetry in the precision attained.

\begin{figure}
\includegraphics[bb=156 194 490 822]{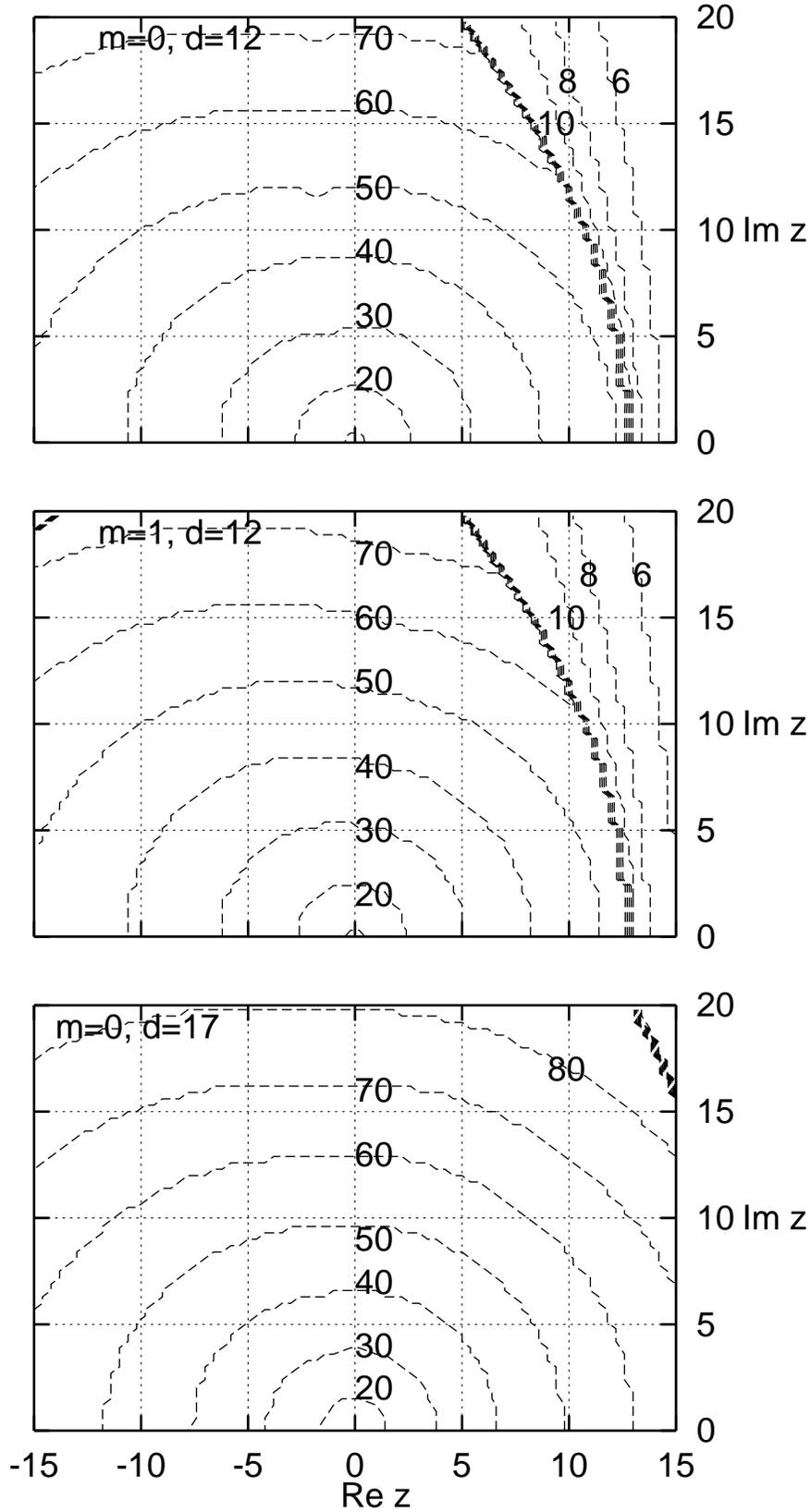}
\caption{
Contour levels of the number $n$ of terms needed to obtain
an accuracy of $d=12$ or $17$ digits of $F_m(z)$, if the power series
(\ref{eq.KummTay}) and the Laurent series (\ref{eq.KummLau})
are used complimentarily.
\label{fig.referenc}
}
\end{figure}

The contrasting regions of good convergence manifested in
Fig.\ \ref{fig.lauren} and Fig.\ \ref{fig.tayl} suggest to combine
these results into Fig.\ \ref{fig.referenc}.
The maximum number of terms needed this way
to obtain $d=12$ digits for $m=0$ 
is $n\le 92$ for the entire $z$ plane (Fig.\ \ref{fig.referenc} top),
to obtain $d=12$ digits for $m=1$ 
is $n\le 84$ for the entire $z$ plane (Fig.\ \ref{fig.referenc} middle),
and to obtain $d=17$ digits for $m=0$ 
is $n\le 166$ for the entire $z$ plane (Fig.\ \ref{fig.referenc} bottom).

\clearpage

\subsection{Gargantini-Henrici Converging Factors}
\begin{figure}
\includegraphics[bb=156 398 561 824]{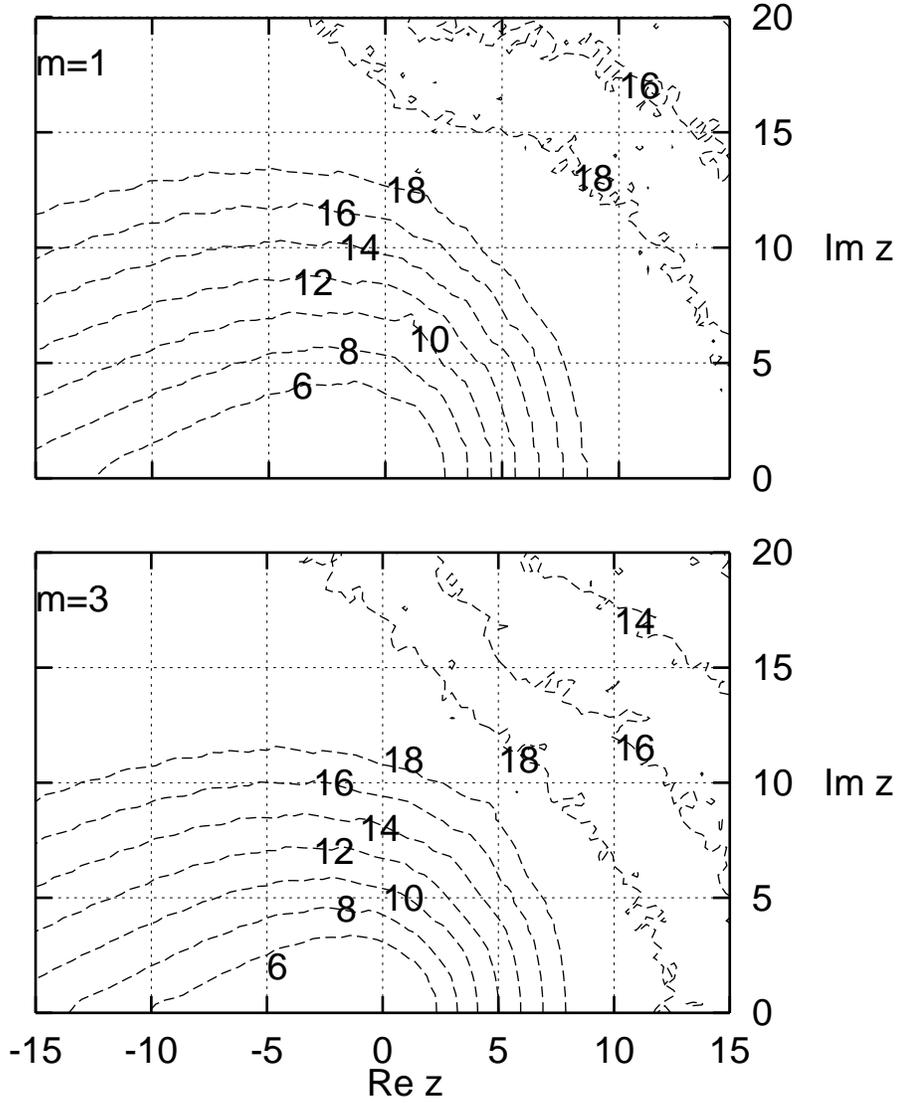}
\caption{
Contour levels of the number of valid decimal digits 
$d=6\ldots 18$ of $F_m(z)$, if the Laurent series (\ref{eq.KummLau})
is summed up to $N=|z|+a$ and the convergent factor $\vartheta_N$
terminated as described in the text.
\label{fig.gargant}
}
\end{figure}

The Gargantini-Henrici analysis \cite{Gargantini}
of the converging factor of the Laurent series (\ref{eq.KummLau})
allows a more accurate calculation of the truncated series
\begin{equation}
\frac{1}{z}\sum_{n=0}^\infty (1-a)_n z^{-n}\equiv\sum_{n=0}^\infty \frac{c_n}{z^{n+1}}.
\label{eq.Stielt}
\end{equation}
The coefficients $c_n\equiv (1-a)_n$ are fed into the quotient-difference
scheme easily derived from \cite[Sec.\ 5]{Gargantini} or taken
from \cite[(3.9.3)]{Wynn1960},
\begin{eqnarray}
q_1^{(n)}&\equiv&\frac{c_{n+1}}{c_n}=1-a+n; \qquad (n=0,1,2,\ldots) \\
e_k^{(n)}&=&k ; \qquad (n=0,1,2,\ldots) \\
q_k^{(n)}&=&n+k-a ;\qquad (n=0,1,2,\ldots)
\end{eqnarray}
and (\ref{eq.Stielt}) is approximated by
\begin{equation}
\sum_{n=0}^{N-1}\frac{c_n}{z^{n+1}}+\frac{c_N}{z^{N+1}}\vartheta_N(z).
\label{eq.Stielt2}
\end{equation}
\begin{equation}
\vartheta_N(z)=\frac{z}{z-}\;\frac{q_1^{(N)}}{1-}\;
\frac{e_1^{(N)}}{z-}\;
\frac{q_2^{(N)}}{1-}\;
\frac{e_2^{(N)}}{z-}\cdots\qquad (N=0,1,2,\ldots)
\label{eq.Garg}
\end{equation}
The following results of Fig.\ \ref{fig.gargant} are based on a
``best knowledge'' approach 
in the sense that the approximation (\ref{eq.Stielt2}) sums to
the same $N\le |z|+a$ as in the previous section, and that the
continued fraction are accumulated until the lower indexes
in $q_k^{(N)}$ and $e_k^{(N)}$ have reached $N$---so to recycle the
same $c_n$ that appear in the main series.
(This roughly triples the number of multiplications and additions
for a particular $z$ compared to the approach of just truncating
(\ref{eq.KummLau}).)

Compared with Fig.\ \ref{fig.lauren}, the multiplication with
the convergent factor has approximately doubled the number of valid digits
in a range of intermediate $|z|$, but as the rational function introduced
by the continued fraction has been allowed to grow to polynomial degrees
of numerator and denominator comparable to the cut-off order of the series,
there is no longer a monotonic increase of accuracy away from the origin.
The regions in the complex plane of predictable accuracy have got a
complicated shape.

\section{Nonstandard Power series}
\subsection{Power Series of the Square}
The square of the series (\ref{eq.Kumm}) may be written with 
the Cauchy product formula \cite[(0.316)]{GR}
as
\begin{equation}
F_m^2(z)=\sum_{n=0}^\infty (-z)^n
\sum_{k=0}^n\frac{1}{k!(n-k)!(2m+2k+1)(2m+2n-2k+1)}
\label{eq.Fmsq1}
\end{equation}
This representation gets slightly more compact by
decomposition into partial fractions and use of the symmetry
$n\leftrightarrow n-k$ in the $k$-sum:
\begin{equation}
F_m^2(z) = 
\sum_{n=0}^\infty \frac{(-z)^n}{2m+n+1}\sum_{k=0}^n \frac{1}{k!(n-k)!(2m+2k+1)}
\label{eq.Fmsq2}
\end{equation}

The decrease of the coefficients of this power series is demonstrated in
Tab.\ \ref{tab.F0sq}. The competitive power series (\ref{eq.KummTay})
converges faster since its coefficients decrease faster, roughly
$\propto 1/(n\cdot n!)$ as a function of $n$.
 
\begin{table}
\begin{tabular}{ccc}
$n$ & $\frac{1}{n+1}\sum_{k=0}^n 1/[k!(n-k)!(2k+1)]$  & $\frac{1}{n+5}\sum_{k=0}^n 1/[k!(n-k)!(2k+5)]$
\\
\hline
0& 1. & (-1)0.40000 00000 00000 00000 00000 00 \\
1& (0)0.66666 66666 66666 66666 66666 67  & (-1)0.57142 85714 28571 42857 14285 71 \\
2& (0)0.31111 11111 11111 11111 11111 11  &  (-1)0.42630 38548 75283 44671 20181 41 \\
3& (0)0.11428 57142 85714 28571 42857 14  & (-1)0.21933 62193 36219 33621 93362 19 \\
4& (-1)0.35132 27513 22751 32275 13227 51  & (-2)0.86974 75364 14203 08086 97475 36 \\
5& (-2)0.93634 76030 14269 68093 63476 03  & (-2)0.28214 74821 47482 14748 21474 82 \\
6& (-2)0.22115 45068 68792 58307 82973 64  & (-3)0.77710 36648 04199 56355 78523 28 \\
7& (-3)0.46965 38029 87136 32046 96538 03  & (-3)0.18637 32720 28586 78606 87262 13 \\
8& (-4)0.90612 20390 19642 50548 12854 38  & (-4)0.39641 62958 93495 68029 29903 10 \\
9& (-4)0.16010 28639 21233 37427 77499 23  & (-5)0.75826 61843 45082 72365 61209 78 \\
10& (-5)0.26077 60505 17784 53196 43665 78  & (-5)0.13186 70614 58730 31544 15735 53 \\
11& (-6)0.39375 88735 04700 13436 03890 04  & (-6)0.21033 98990 99545 68213 83768 14 \\
12& (-7)0.55388 71577 33648 50973 95012 03  & (-7)0.30997 96436 79842 49838 36209 14 \\
13& (-8)0.72901 15029 68230 89769 97160 49  & (-8)0.42463 93655 23994 83198 68341 66 \\
14& (-9)0.90128 71648 37490 71509 80843 99  & (-9)0.54355 16227 46159 40865 55127 90 \\
15& (-9)0.10503 39537 53128 82274 43051 64  &  (-10)0.65304 55397 25446 17367 34616 47\\
16& (-10)0.11574 62528 45873 78175 30636 17  & (-11)0.73931 27958 11477 67351 64672 14 \\
17& (-11)0.12095 88230 02785 23102 13378 56  & (-12)0.79139 24236 39157 01672 60238 06 \\
18& (-12)0.12018 49078 18966 37363 41310 02  & (-13)0.80345 79062 12657 33014 90722 55 \\
19& (-13)0.11380 77808 88498 70226 59025 07  & (-14)0.77576 42229 41899 25361 71656 84 \\
20& (-14)0.10293 02165 82138 49790 10396 95 & (-15)0.71410 21050 93626 04408 94055 80 
\end{tabular}
\caption{
Decrease of the coefficients in (\ref{eq.Fmsq2}), cases $m=0$ and $m=2$,
as a function of $n$. The numbers in parentheses denote multiplication by
powers of 10, e.g., for $n=17$, we have $0.1209\ldots\cdot 10^{-11}$
and $0.791\ldots\cdot 10^{-12}$.
\label{tab.F0sq}
}
\end{table}
 
\subsection{Power Series with Split-Off Exponential}\label{sec.Ftayl2}
The fundamental power series (\ref{eq.Kumm}) is derived from replacing
$e^{-zt^2}$ by its power series; by any truncation of the series,
this represents the integral kernel just at $t=0$.
We investigate a more accurate interpolation, which equals the kernel
at {\em both\/} limits of the $t$-interval $[0,1]$, which
splits off a simpler exponential that still can be integrated exactly,
and which accumulates the (smaller) remainder in 
a different power series:
\begin{equation}
e^{-zt^2}=e^{-zt}+z(t-t^2)-\frac{z^2}{2}(t^2-t^4)+\frac{z^3}{3!}(t^3-t^6)-\cdots .
\end{equation}
\begin{equation}
F_m(z)=\int_0^1t^{2m}e^{-zt}dt
-\sum_{n=1}^\infty \frac{(-z)^n}{(n-1)!}\frac{1}{(2m+n+1)(2m+2n+1)} .
\label{eq.FTayl2}
\end{equation}
In particular,
\begin{equation}
F_0(z)=\frac{1-e^{-z}}{z}
-\sum_{n=1}^\infty \frac{(-z)^n}{(n-1)!}\frac{1}{(n+1)(2n+1)}.
\label{eq.FTayl20}
\end{equation}
As a side note, a decomposition in partial fractions and insertion of the
hypergeometric notation
with (\ref{eq.Kumm}) at $m=1$ yields
\begin{eqnarray}
F_0(z)
&=&\frac{1-e^{-z}}{z}
+z\sum_{n=0}^\infty \frac{(-z)^n}{n!}\frac{1}{(n+2)(2n+3)} \\
&=&\frac{1-e^{-z}}{z}
+z
\underbrace { \sum_{n=0}^\infty \frac{(-z)^n}{n!}\frac{1}{n+3/2} }_{\frac{2}{3}{_1}F_1(3/2;5/2;-z)=2F_1(z)}
-z
\underbrace { \sum_{n=0}^\infty \frac{(-z)^n}{n!} \frac{1}{n+2} }_{\frac{1}{2}{_1}F_1(2;3;-z)=[1-(z+1)e^{-z}]/z^2},
\end{eqnarray}
and emerges as a complicated adornment of (\ref{eq.recurr}).
Supposed one has a fast, reliable method to compute $1-e^{-z}$,
(\ref{eq.FTayl2}) looks beneficial compared to (\ref{eq.KummTay}) because
the total power of the summation variable $n$ in the denominator is 
slightly larger. Graphing the results of (\ref{eq.FTayl20}) the same way as in the two
upper plots of Fig.\ \ref{fig.tayl}, however, would yield no differences
visible to the eye.
Therefore we do not look into this ansatz further.

\subsection{Power Series of the Half Argument}
Convergence of power series is generally faster closer to the origin;
the trigonometric identity $\sin^2(\omega/2)=[1-\cos\omega]/2$ allows us to
reduce the distance between $z$ and the origin by half
if we substitute $t=\sin(\omega/2)$ in (\ref{eq.FmInt}):
\begin{equation}
F_0(z)=\frac{1}{2}\int_{-1}^1 e^{-zt^2}dt=\frac{e^{-z/2}}{4}\int_{-\pi}^\pi
e^{\frac{z}{2}\cos \omega}\cos{\frac{\omega}{2}}d\omega
=
\frac{e^{-z/2}}{4}\sum_{n=0}^\infty \frac{1}{n!}\left(\frac{z}{2}\right)^n
\int_{-\pi}^\pi \cos^n\omega \cos\frac{\omega}{2}d\omega.
\end{equation}
The auxiliary $\omega$-integrals would be drawn from the recursion \cite[(2.538.1)]{GR}
\begin{equation}
\int_{-\pi}^\pi \cos^n\omega \cos\frac{\omega}{2}d\omega
=\frac{1}{n+1/2}\left[ 2(-)^n+n
\int_{-\pi}^\pi \cos^{n-1}\omega \cos\frac{\omega}{2}d\omega
\right].
\end{equation}
Since these are of the order of 1 for all $n$, we are left with a power series
which converges $\sim z^n/(n!2^n)$, which is to be compared to
$\sim z^n/(n!n)$ of (\ref{eq.KummTay}).

\begin{figure}
\includegraphics[bb=144 197 496 820]{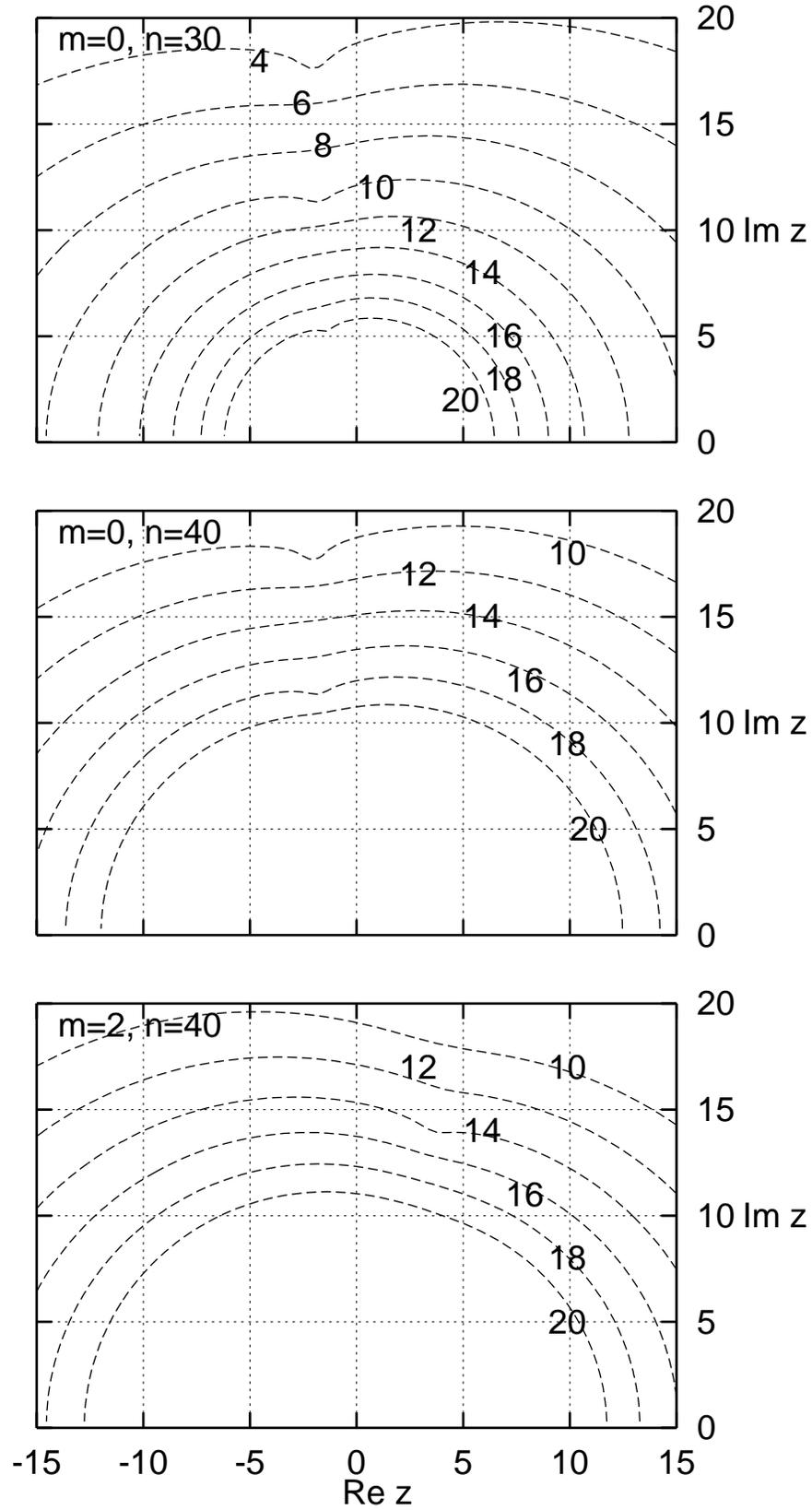}
\caption{
Contour levels of the number of valid decimal digits 
$d=4\ldots 20$ of $F_m(z)$ by the power series (\ref{eq.taylH}),
if the power series is truncated after $n=30$ or $40$ terms.
For $m=0$ and $n=60$, $d$ is  $\ge19.6$ in this $z$-domain.
\label{fig.taylHalf}
}
\end{figure}

The generalization to nonzero $m$ reads
\begin{equation}
F_m(z)
=
\frac{e^{-z/2}}{4}\sum_{n=0}^\infty \frac{1}{n!}\left(\frac{z}{2}\right)^n
\int_{-\pi}^\pi \sin^{2m}\frac{\omega}{2}\cos^n\omega \cos\frac{\omega}{2}d\omega .
\label{eq.taylH}
\end{equation}
The ($n,m$)-table of the $\omega$-integrals could be generated from the
table at $m=0$ or from scratch,
\begin{eqnarray}
\int_{-\pi}^\pi \sin^{2m}\frac{\omega}{2}\cos^n\omega \cos\frac{\omega}{2}d\omega
&=&
\frac{1}{2^m}\sum_{k=0}^m (-)^k{m \choose k}\int_{-\pi}^\pi \cos^{k+n}\omega \cos\frac{\omega}{2}d\omega \nonumber\\
&=&
\sum_{k=0}^m\sum_{l=0}^n{m \choose k}{n \choose l}(-)^{k+n+l}2^{l+1}\frac{(k+l)!}{(\frac{1}{2})_{k+l+1}}
\label{eq.omegnm}
\end{eqnarray}
The examples of Fig.\ \ref{fig.taylHalf} demonstrate that
about 8 additional digits have been gained for the case
$m=0$ and $n=30$ terms relative to the data of Fig.\ \ref{fig.tayl}.
There is no additional run-time cost since the coefficients table
of (\ref{eq.omegnm}) is static without $z$-dependence.

\clearpage
\section{Interpolating the Index}

The integral (\ref{eq.FmInt}) is solvable if $2m=2a-1$
is an {\em odd\/} integer \cite[(2.321.2)]{GR}:
\begin{equation}
F_m(z)=\frac{(a-1)!}{2z^a}-\frac{e^{-z}}{2z}\left[ 1+\sum_{n=1}^{a-1}
\frac{(1-a)(2-a)\cdots(n-a)}{(-z)^n}\right]
, (a=1,2,3,\ldots).
\label{eq.Fofm}
\end{equation}
This could also be derived from (\ref{eq.Kumm}) and (\ref{eq.KummLau}), where
the sum in
(\ref{eq.KummLau}) terminates if $a$ is a positive integer.

Given $z$, let (\ref{eq.Fofm}) be computed for an index set
$m=1/2, 3/2, 5/2,\ldots, N+1/2$, optionally facilitated by
the index recursion \cite[(2.321.1)]{GR} or (\ref{eq.recurr}).
Let us pursue the idea that approximate values at intermediate $m=1, 2,\ldots$,
which we are actually interested in, are deduced by some interpolation.

\begin{figure}
\includegraphics[bb=116 218 377 761]{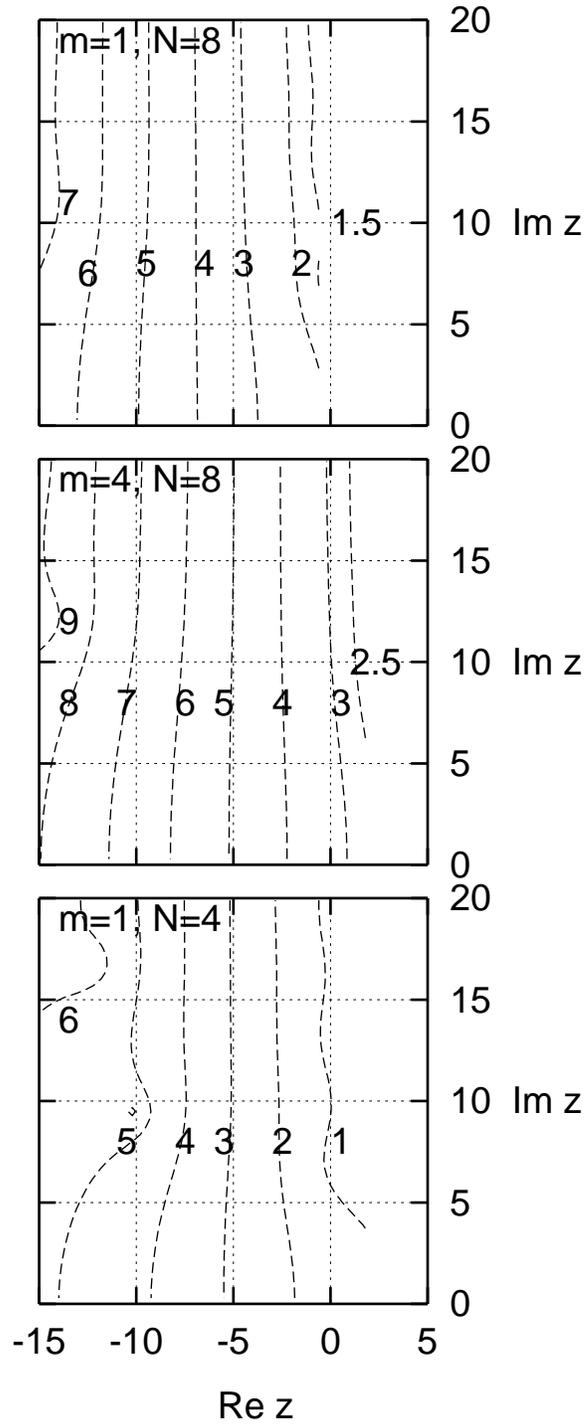}
\caption{
Contour levels of the number of valid decimal digits 
$d=1\ldots 9$ of $F_1(z)$ and $F_4(z)$ by
interpolation of the index with a polynomial of degree $N=8$ or $4$ in $m$,
Eq.\ (\ref{eq.polym}).
\label{fig.Fofm}
}
\end{figure}

If this is done by the unique interpolating polynomial
$\sum_{j=0}^Nb_jm^j$
of degree $N$ in $m$, the intermediate step of the calculation
can formally be written down as a $(N+1)\times(N+1)$ inhomogeneous
system of linear equations to get the $N+1$ unknown polynomial coefficients
$b_j$ (a Lagrange Interpolation might be cheaper numerically, though),
\begin{equation}
\sum_{j=0}^N m_k^j b_j = F_{m_k}(z),\qquad (k=0,1,\ldots,N; m_k=k+1/2),
\label{eq.polym}
\end{equation}
and we get for
example Fig.\ \ref{fig.Fofm}.
We see: (i) The method becomes more and
more unreliable as $\Re z$ increases, which emerges from different weighting
by the factors $t^{2m}$ and $e^{-zt^2}$ in the integral kernel:
The derivative with respect to $m$ multiplies the integrand (\ref{eq.FmInt})
by $2\ln t$. To keep the derivative small relative to $F_m(z)$ itself
(and to keep $F_m(z)$ a flat function of $m$), it is
advantageous that $t^{2m}e^{-zt^2}$ weights stronger  at the right
limit of the $t$-interval, where $\ln t$ stays small.
$t^{2m}$ becomes large near $t=1$, whereas $e^{-zt^2}$ is larger near $t=0$
or $t=1$ depending on the sign of $\Re z$.)
(ii) The estimates are better close to the middle of the $m$
interval $[\frac{1}{2},N+\frac{1}{2}]$, which
was sampled to define the interpolation polynomial, than for
$m$ close to the interval limits---which is expected for any interpolation derived from
approximately equidistant sampling points.

No such interpolation polynomial of $m$ satisfies (\ref{eq.recurr})
on the entire real $m$-axis---though the exact $F_m(z)$ does for all $m>1/2$.
Therefore one could seek after improvement of this interpolation
by enforcing some compliance with (\ref{eq.recurr}).
Test calculations were made following the strategy that one or more
lines in the system of linear equations (\ref{eq.polym})
are replaced by coupling some of the unknown $F_m(z)$ via (\ref{eq.recurr}),
in concrete
\begin{equation}
\sum_{j=0}^N [2zm_u^j-(2m_u-1)m_v^j] b_j = -e^{-z}, \qquad
(m_v=0,1,2,\ldots; m_u=m_v+1).
\label{eq.mrecurr}
\end{equation}
The polynomial that ensues does no longer hit the $F_{m_k}(z)$ that was 
thrown out, but follows (\ref{eq.recurr}) for those $m_u$ brought in.
Results of this ansatz with $N=4$, removal of $F_{9/2}$ from the set
of interpolated points, and introduction of the $m_v=1\leftrightarrow
m_u=2$ coupling with (\ref{eq.mrecurr}) add at most about half a digit of
accuracy to what is already shown in Fig.\ \ref{fig.Fofm}.
There is no further improvement if the line for $F_{7/2}$ in
(\ref{eq.polym}) is also removed to add the line (\ref{eq.mrecurr})
for the $m_v=2\leftrightarrow m_u=3$ coupling.

Furthermore one would try to add known information on the first derivatives
with respect to $m$ to enhance the quality of the polynomial interpolation.
Unfortunately, even the simplest case of $\partial F_m(z)/\partial m$ at
$m=1/2$ would demand computation of the Exponential Integral $E_1(z)$ \cite{Cody},
\begin{equation}
\frac{\partial F_m(z)}{\partial m}_{|m=1/2}
=-\frac{1}{2z}\left[ E_1(z)+\gamma+\log(z)\right].
\end{equation}
This infects also the partial derivatives at $m=3/2, 5/2,\ldots$
\begin{equation}
\frac{\partial}{\partial m}F_m(z)=\frac{1}{z}F_{m-1}(z)+\frac{m-1/2}{z}
\frac{\partial}{\partial m}F_{m-1}(z); \qquad (m>1/2).
\end{equation}

Another attempt of refinement is to acknowledge the 
simple poles at the negative half integers of $m$. Test calculations
with the modified separation ansatz \cite{Berrut}
$F_m(z)\approx (\sum_{j=0}^Nb_jm^j)/(2m+1)$, which manifests the
pole at $m=-1/2$,
and again a numerator polynomial of degree $N=4$ result in
changes of up to one digit (in both directions)
compared to the bottom graph in Fig.\ \ref{fig.Fofm}.

In summary, it seems to be difficult to bridge the gap between
$F_m$ at half integer and integer $m$ through interpolation.

\clearpage

\section{Generic Methods of Integration}
\subsection{Local Taylor Expansions in the Integration Interval}\label{herm.sec}
\subsubsection{Expansion of the Exponential}
\begin{figure}
\includegraphics[bb=184 200 561 824]{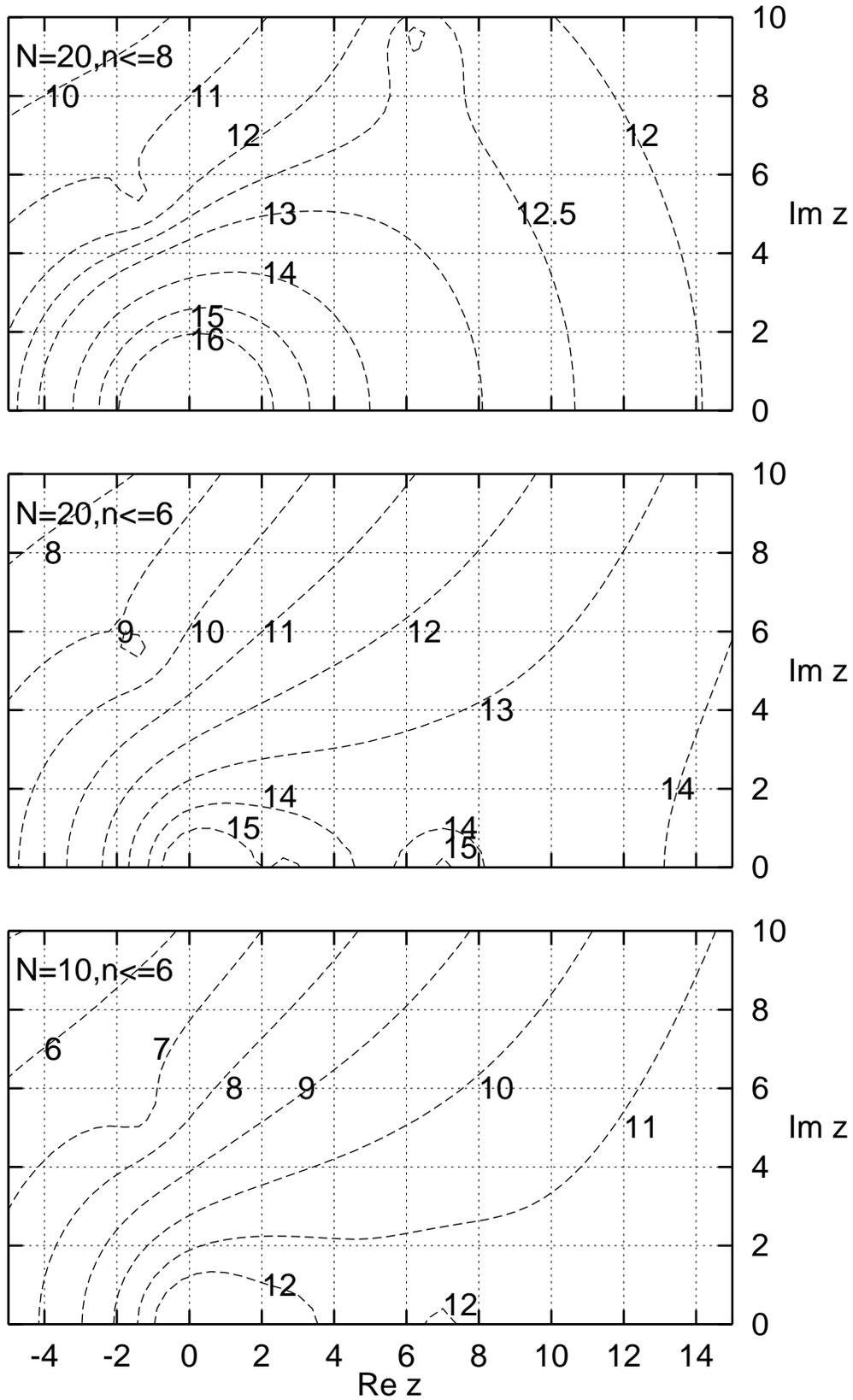}
\caption{
Contour levels of the number of valid decimal digits 
$d=8\ldots 16$ of $F_0(z)$, if the series (\ref{eq.herm})
is used with $N=10$ or $20$ subintervals and if the Hermite Polynomials
are included up to degree $n=6$ or $8$. Note that in a region
of $\Re z\gtrapprox 8$ and small $|\Im z|$ an increase of $n$ from $6$ to $8$
decreases the accuracy.
\label{fig.Herm}
}
\end{figure}

The ``global'' first approximation of the integrand in Sec.\ \ref{sec.Ftayl2}
may be pushed one notch further towards a brute-force numerical method
by slicing the $t$-interval $[0,1]$ into $N$ same size subintervals
of half-width $\Delta\equiv 1/(2N)$, centered at $t_l=(2l-1)/\Delta$
($l=1,\ldots,N$). In each of these subintervals,
$\exp(-zt^2)$
is approximated by its Taylor series around $t_l$.
The derivatives are \cite[(7.1.19)]{AS}
\begin{equation}
\left(\frac{d}{dt}\right)^ne^{-zt^2}=(-)^nz^{n/2}H_n(\sqrt{z}t)e^{-zt^2},
\end{equation}
in terms of Hermite Polynomials $H_n$---looking at the exponential as the
first derivative of the error function---,
whence the Taylor series
\begin{equation}
e^{-zt^2}=\sum_{n=0}^\infty \frac{1}{n!}(t-t_l)^n(-)^n z^{n/2}H_n(\sqrt{z}t_l)
e^{-zt_l^2}.
\end{equation}
$F_0(z)$ is the Riemann sum over the subintervals
\begin{equation}
F_0(z)\approx\sum_{l=1}^N
\int_{t_l-\Delta}^{t_l+\Delta}dt e^{-zt^2}=2
\sum_{l=1}^N e^{-zt_l^2}
\sum_{n=0,2,4,\ldots}^\infty \frac{\Delta^{n+1}}{(n+1)!} z^{n/2}H_n(\sqrt{z}t_l)
\label{eq.herm}
\end{equation}
Fig.\ \ref{fig.Herm} shows the accuracy of (\ref{eq.herm}) for two 
$N$, keeping the polynomial expansion only up to some degree $n$.

The two graphs in Fig.\ \ref{fig.Herm} with degrees kept up to $n=6$
show the two spots at $z\approx 2.80$ and $z\approx 7.02$ with higher precision
than their surroundings. This is part of a more general phenomenon, in which
for some $z$ the positive and negative lobes of the Hermite Polynomials are
sampled with best effective cancellation. This generates one such spot
at $z\approx 4.08$ ($n\le 4$), those two for $n\le 6$,\ldots

On the computational expense: as the formula requires Hermite polynomials of even indexes only,
there is no need to compute the $\sqrt z$ and no gain using their
recursion formulas. Two (complex) multiplications
compute $zt_l^2$. Each $H_n$ is a polynomial of degree $n/2$ in this combined
variable, which costs  $n/2$ multiplications and $n/2$ additions with the Horner scheme.
The powers of $\Delta$ are fixed, and there are about 3 multiplications
for each term in $n$, plus one (complex) exponential. This is to be
multiplied by $N$.
This totals at least $16N$ multiplications if
$n$ runs up to 6, and so only the lowermost picture in Fig.\ \ref{fig.Herm}
would be part of a fair comparison with the other competitive 
approaches.

\subsubsection{Expansion of the Algebraic Factor}
In a similar manner as above, one could expand the algebraic factor $u^{m-1/2}$
of (\ref{eq.FmInt}) in Taylor series around $u_l\equiv (2l-1)/\Delta$ to end
up with closed form integrals. To second order,
\begin{equation}
u^{m-1/2}\approx u_l^{m-1/2}+(m-\frac{1}{2})u_l^{m-3/2}(u-u_l)
+(m-\frac{1}{2})(m-\frac{3}{2})u_l^{m-5/2}\frac{(u-u_l)^2}{2},
\label{eq.utayl}
\end{equation}
and the Riemann sum
\begin{eqnarray}
F_m(z)\approx-\frac{1}{2z}\sum_{l=1}^N &&
e^{-zu_l}u_l^{m-1/2}\bigg\{ 
e^{-z\Delta}-e^{z\Delta}+\frac{m-1/2}{zu_l}
\left[e^{-z\Delta}(1+z\Delta)-e^{z\Delta}(1-z\Delta)\right] \nonumber \\
&& +\frac{(m-1/2)(m-3/2)}{(zu_l)^2}
\left[e^{-z\Delta}(\frac{z^2\Delta^2}{2}+z\Delta+1)
-e^{z\Delta}(\frac{z^2\Delta^2}{2}-z\Delta+1)\right]
\bigg\} . \nonumber
\end{eqnarray}
This achieves up to
$d=2.6$ decimals ($m=0$, $N=20$),
$d=5.4$ ($m=1$, $N=20$),
$d=5.8$ ($m=1$, $N=40$), and $d=8.2$ ($m=2$, $N=40$)
in the $z$-domain as in Fig.\ \ref{fig.Herm}.
(These numbers
refer to the ``lower left'' corner of the $z$-region, and are a few digits
worse in the opposite corner.)

The main obstacle to higher performance is the poor fit of (\ref{eq.utayl})
close to $u=0$. We may patch this by replacing the contribution in this subinterval,
the term $l=1$ where $0\le u \le 2\Delta$ remains small, by the associated
power series of the exponential,
\begin{equation}
F_m(z)\approx-\frac{1}{2z}\sum_{l=2}^N 
e^{-zu_l}u_l^{m-1/2}\bigg\{ \ldots
\bigg\} +\frac{(2\Delta)^{m+1/2}}{2}\sum_{n=0,1,2,\ldots}\frac{(-2z\Delta)^n}{(m+1/2+n)n!}.
\end{equation}
With this ansatz and the sum over $n$ kept up to $n=3$, the maximum
number of digits in the $z$-domain as in Fig.\ \ref{fig.Herm} rise to
$d=7.8$ ($m=1$, $N=40$),
and $d=10.4$ ($m=2$, $N=40$).

\subsection{Fourier Expansion of the Algebraic Integral Kernel}
\begin{figure}
\includegraphics[bb=184 200 561 824]{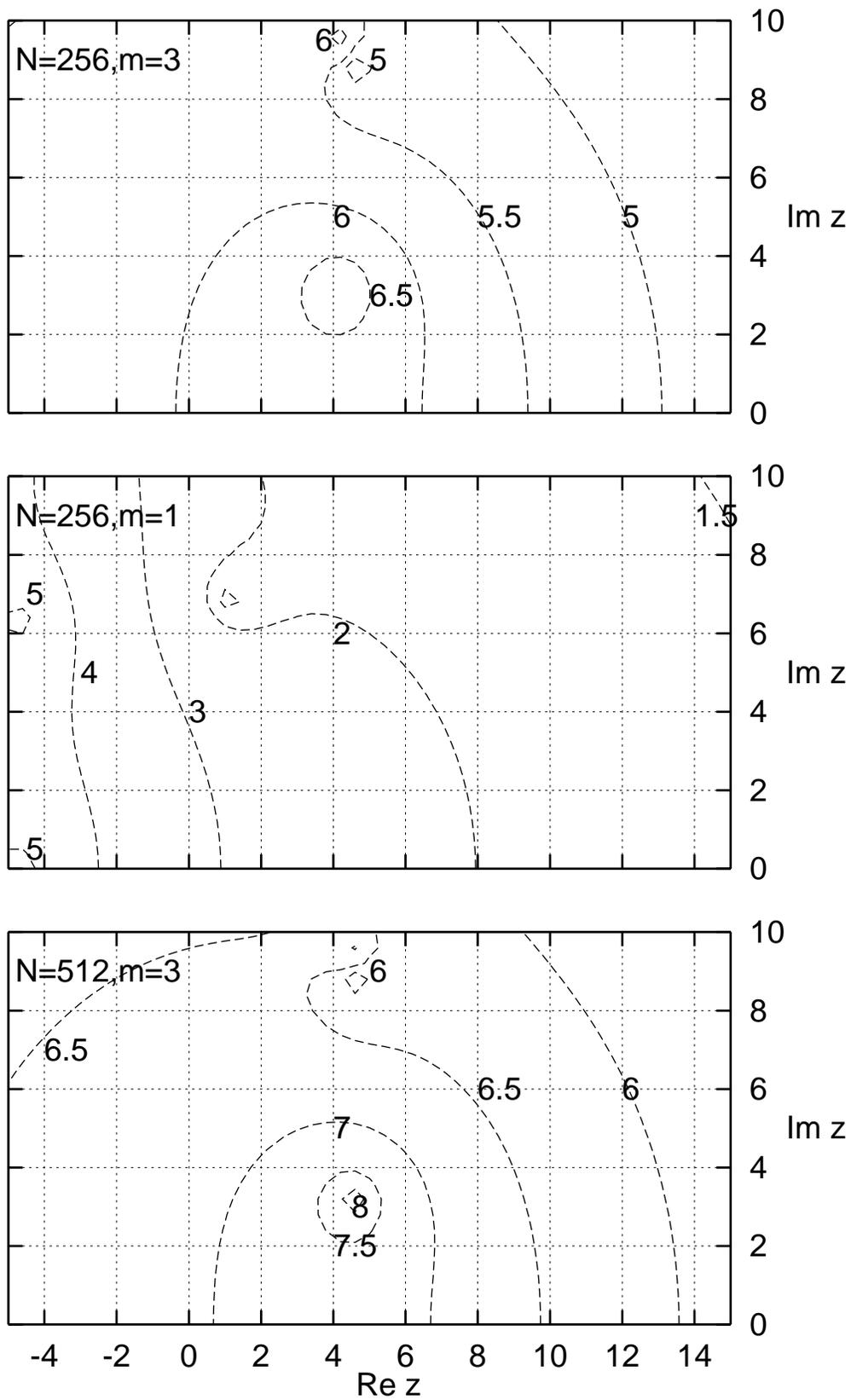}
\caption{
Contour levels of the number of valid decimal digits 
$d=1.5\ldots 8$ of $F_m(z)$, if the Fourier series (\ref{eq.uFour2})
is used with $N=256$ or $512$ abscissa points of the Fourier interpolation.
\label{fig.ft}
}
\end{figure}

A Fourier expansion of the algebraic term in (\ref{eq.FmInt})
\begin{equation}
u^{m-1/2}=\sum_{l=0}^{\infty}c_l \cos(lu\frac{\pi}{2})
\label{eq.uFour}
\end{equation}
offers the series
\begin{equation}
\int_0^1u^{m-1/2}e^{-zu}du=\frac{1-e^{-z}}{z}+2\sum_{l=1,3,5,\ldots}^\infty c_l
\frac{2z+(-)^{[l/2]}l\pi e^{-z}}{(2z)^2+(l\pi)^2},
\label{eq.uFour2}
\end{equation}
where $c_0=1$ and $c_2=c_4=c_6=\ldots=0$ have already been assumed.
To reduce any Gibbs oscillations \cite{Driscoll} of (\ref{eq.uFour}) at the
ends
of the interval $[0,1]$, $u^{m-1/2}$ is embedded into the even, 4-periodic, and steady
carrier function
$f(u)\equiv u^{m-1/2}$ for $u\in [0,1]$, $f(u)\equiv 2-(2-u)^{m-1/2}$ for $u\in [1,2]$,
$f(u)\equiv f(2-u)$ and $f(u)=f(u+4)$ elsewhere.
(Obviously, the singularity at $u=0$ reduces the quality of this approach
right from the start if $m=0$.)
The $c_l$ are
approximated by a discrete cosine transform on $N$ grid points
\begin{equation}
c_j=\frac{4}{N}\left[ (-)^j +\sum_{k=1}^{N/2-1}f(4k/N)\cos(2\pi\frac{jk}{N})\right]\qquad (j=1,2,\ldots,N/2),
\end{equation}
and the summation (\ref{eq.uFour2}) is truncated at $l=N/2$.
The $c_l$ would be kept in constant tables since they do not depend on $z$,
and the cost of evaluating (\ref{eq.uFour2}) amounts to about $N/4$ evaluations
of the rational term.
Fig.\ \ref{fig.ft} shows that for moderately small $m$ a precision of just of the order
of 5 digits result from this type of evaluation, which is not efficient
compared to other methods proposed here
This is ultimately a progression of residual fitting
errors (Gibbs oscillations) which remain rather large
close to $u=0$ and $u=1$ (Tab.\ \ref{tab.ft}).

\begin{table}
\begin{tabular}{cccc}
$m$ & \multicolumn{3}{c}{$N$} \\
 & 512 &          256              & 128 \\
\hline
1  & $(-2)3.6$  & $(-2)5.1$ & $(-2)7.4$ \\
2  & $(-5)4.6$  & $(-4)1.3$ & $(-4)3.7$ \\
3  & $(-6)8.2$  & $(-5)3.5$ & $(-4)1.4$
\end{tabular}
\caption{
Typical maximum deviation between the left hand side and the right
hand side of (\ref{eq.uFour}) in the interval $u\in [0,1]$,
if one period of $f(u)$ is sampled by $N$ points
and the sum (\ref{eq.uFour})  limited to $0\le l\le N/2$.
Numbers in parentheses are powers of $10$ as in Tab.\ \ref{tab.F0sq}.
\label{tab.ft}
}
\end{table}

In the special case where one would like to tabulate $F_m(z)$ along some
fixed ordinate $\Re z=const$ in equidistant steps of $\Im z$,
one could draw a lot of additional benefit
from the inherent parallelism of Fast Fourier Transform techniques
applied to the decomposition
\begin{equation}
F_m(z)=\frac{1}{2}\int_0^1 u^{m-1/2}e^{-\Re zu}e^{-i\Im zu}du.
\end{equation}
This deems to be too special to be put into detail here.

\clearpage
\subsection{Taylor series in the complex $z$-Plane}\label{sec.kumm3}

With \cite[(13.4.9)]{AS}, the complex derivatives of $F_m(z)$
are equivalent to ladder-type operation with respect to $m$ and $a$:
\begin{equation}
\left(\frac{d}{dz}\right)^n {}_1F_1(a;a+1;z)=\frac{a}{a+n}{}_1F_1(a+n;a+1+n;z),
\end{equation}
\begin{equation}
\left(\frac{d}{dz}\right)^n F_m(z)=(-)^nF_{m+n}(z)
\label{eq.diff}
\end{equation}
This close link between the index and the higher derivatives means
that one may tabulate the expansion coefficients of the Taylor series
\begin{equation}
F_m(z_0+\varepsilon) = F_m(z_0)-\varepsilon F_{m+1}(z_0)+\frac{\varepsilon^2}{2!}F_{m+2}(z_0)
+\cdots+\frac{(-\varepsilon)^n}{n!}F_{m+n}(z_0)+\cdots
\label{eq.kumm3}
\end{equation}
anchored at some $z_0$ in the complex plane for all $m$ at the same time.
This keeps these tables smaller than for any function with 
decoupled index and argument. Despite the fact that these tables need
to contain complex values and need to be arranged on a 2D grid in the
complex plane, there is no inherently new aspect over what is already
assessed in the literature for real $z$ \cite{Takashima}.
Chebychev approximations in the complex case
could be derived by expanding $(-\varepsilon)^n$
in a sum over products of $\Re \varepsilon$ and $\Im \varepsilon$, normalization
to the interval $[-1,1]$ and independent progression for the real and imaginary
part as described in \cite[\S 22.20]{AS}.

\begin{figure}
\includegraphics[bb=73 93 509 372]{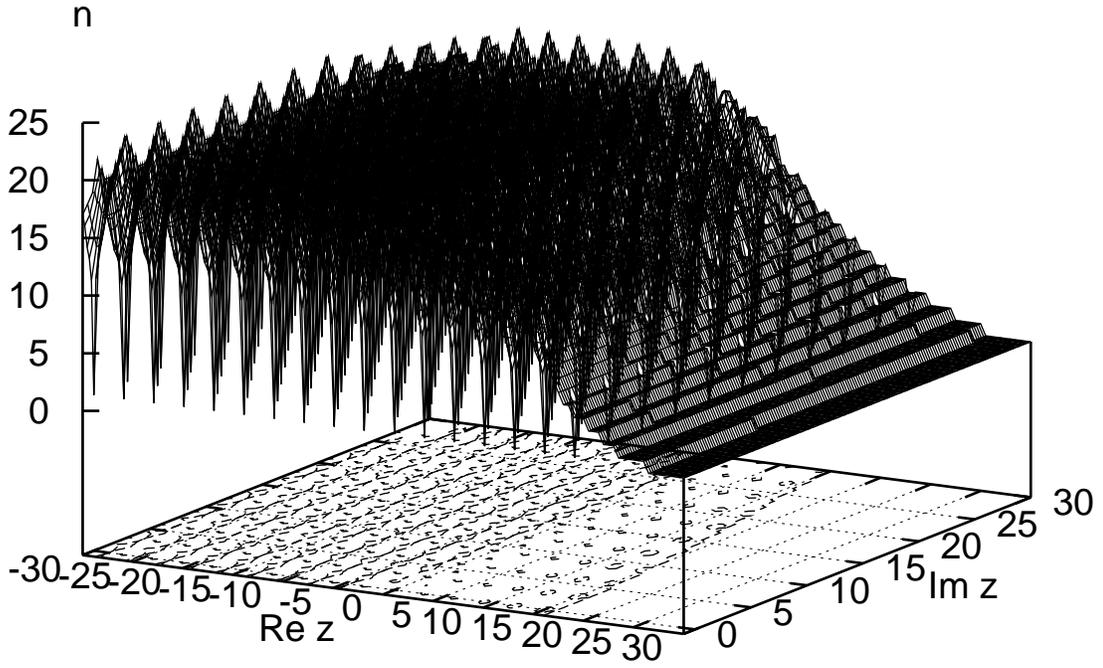}
\caption{
The number $n$ of terms in the Taylor series Eq.\ (\ref{eq.kumm3})
to achieve $d=14$ digits of accuracy of $F_0(z)$.
\label{fig.kumm3T0}
}
\end{figure}

Fig.\ \ref{fig.kumm3T0} is an example where the nodes span the $\Re z_0$
interval from $-33$ to $18$ with a stride of $s=3$, and $\Im z_0$
the interval from $0$ to $36$ also with a stride of $s=3$:
\begin{equation}
z_0= ks+ils\qquad (k=-11,-10,\ldots,6; l=0,1,\ldots,12)
\end{equation}
Up to 23 terms
need to be accumulated in (\ref{eq.kumm3}) to calculate $F_0(z)$ if $z$ falls
inside this finite grid's domain;
values outside are handled with Eq.\ \ref{eq.KummLau}.
This term count is stored as entry $23$ at $m=0$ and $d=14$ in the upper part
of Tab.\ \ref{tab.kumm3}.
A free entry in the table indicates,
what could not be computed if the tabulated $F_j(z)$ would be limited to $j\le 30$,
which would be less than $m+n$ for this slot.
The lower part of the table illustrates, how a three times denser grid
of nodes
cuts down on the worst case convergence, as it reduces the
maximum distance $|\varepsilon|$ to the nearest
$z_0$ from $3/\sqrt 2$ to $1/\sqrt 2$.
This necessitates a ninefold larger static table, unless one uses the 
smaller number of terms to reduce the maximum index $j$ of tabulated $F_j$.

\begin{table}
\begin{tabular}{ccccccc}
$m$ & \multicolumn{5}{c}{$d$} \\
 & 12 & 14 & 15 & 16 & 17 \\
\hline
0  & 21 & 23 & 24 & 25 & 26 \\
1  & 21 & 23 & 24 & 25 & 26 \\
3  & 21 & 23 & 24 & 25 & 26 \\
5  & 21 & 23 & 24 & 25 & \\
\end{tabular}

\begin{tabular}{ccccccc}
$m$ & \multicolumn{5}{c}{$d$} \\
 & 12 & 14 & 15 & 16 & 17 \\
\hline
0  & 14 & 16 & 17 & 17 & 18 \\
1  & 14 & 16 & 16 & 17 & 18 \\
5  & 14 & 15 & 16 & 17 & 18 \\
\end{tabular}
\caption{
The upper table shows the maximum number of terms needed for nearest neighbor
Taylor expansions (\ref{eq.kumm3}) as a function of $m$ and accuracy $d$
on a support grid with spacing 3 in the $z$ plane.
The lower table illustrates the case of a denser spacing of 1.
\label{tab.kumm3}
}
\end{table}

On a side note, (\ref{eq.kumm3}) could be rewritten
similar to \cite[(2.1)]{Gautschi2003}
\begin{equation}
F_m(z_0+\varepsilon)=
e^{-z_0}F_m(\varepsilon)+z_0\sum_{n=0}^{\infty}\frac{1}{m+1/2+n}\frac{(-\varepsilon)^n}{n!}F_{m+1+n}(z_0),
\end{equation}
using (\ref{eq.recurr}) once for each term, or integrating (\ref{eq.FmInt})
by parts.

\clearpage

\subsection{Gauss-Jacobi quadrature}\label{GaussJ.sec}
As considered by Gautschi \cite{Gautschi2002} for a more general case,
(\ref{eq.FmInt}) is readily accessible by a Gauss-Jacobi quadrature,
\begin{equation}
F_m(z) \approx \sum_{i=1}^n w_i e^{-zt_i^2} .
\label{eq.gaussJ}
\end{equation}
Weights
$w_i$ and abscissae $t_i$ are discussed in App.\ \ref{app.gaussJ}.

\begin{figure}
\includegraphics[bb=151 196 496 820]{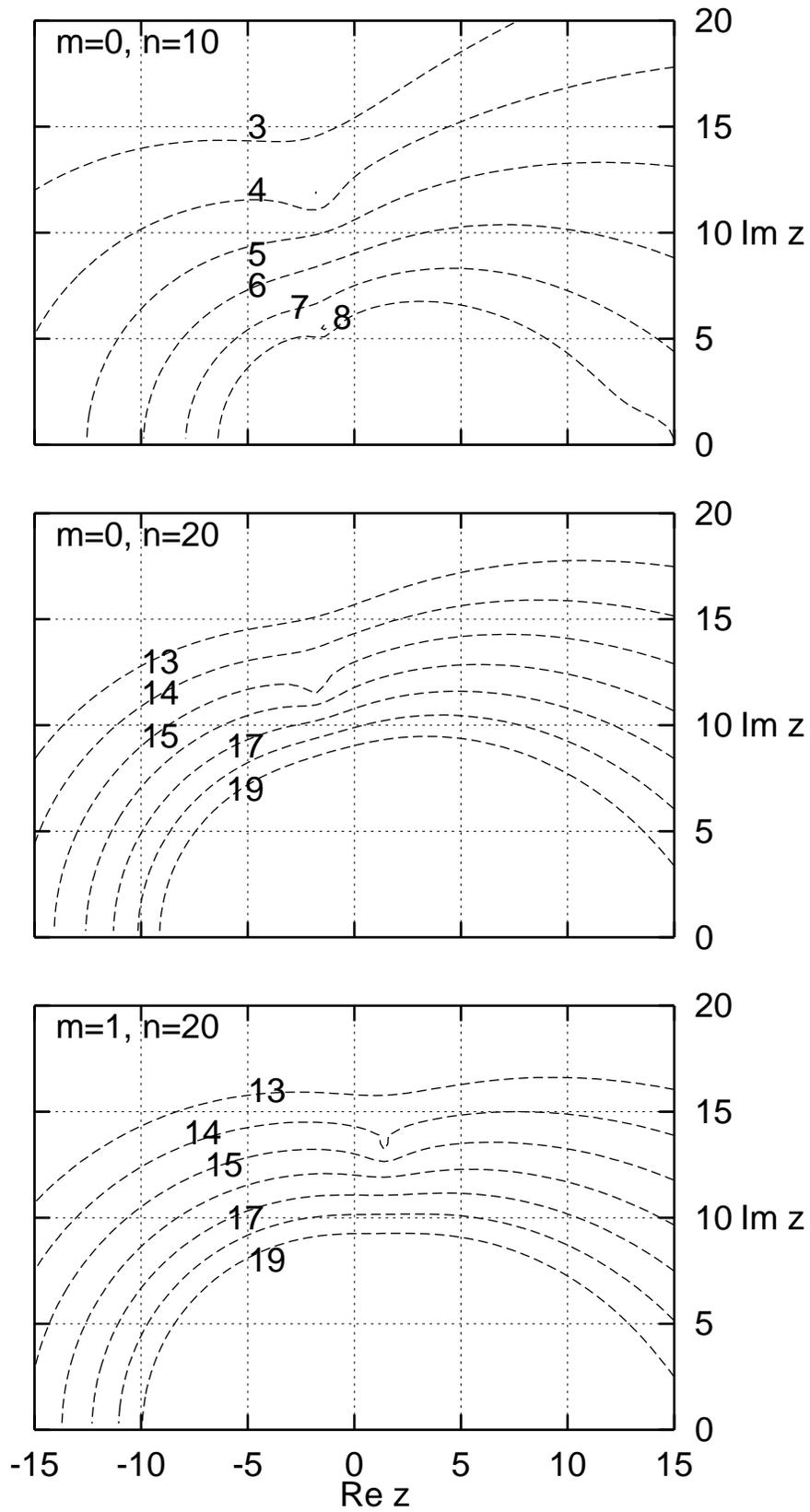}
\caption{
Contour levels of the number of valid decimal digits 
$d=3\ldots 19$ of $F_m(z)$ by (\ref{eq.gaussJ}), if the order of the
Gauss-Jacobi quadrature is $n=10$ or $20$.
\label{fig.gaussJ}
}
\end{figure}

Fig.\ \ref{fig.gaussJ} demonstrates that this method is most robust at
small $|z|$, which is expected since the Gauss quadrature effectively
approximates $\exp(-zt^2)$ by a polynomial of degree $2n$ in $t$.
The numerical expense roughly adds up to $n$ computations of exponentials
$\exp(-zt_i^2)$. This is cheaper than evaluation of (\ref{eq.herm})
at the same $N=n$; by further comparison of Fig.\ \ref{fig.gaussJ} with
Fig.\ \ref{fig.Herm} we conclude that the Gauss-Jacobi ansatz proposed
here
is superior to the method of Sect.\ \ref{herm.sec}.

The drawback by further comparison with the reference calculation
of Sec.\ \ref{sec.referenc}
is that each complex exponential needs much more CPU time than
a complex multiplication.
Tests with the C++ implementation by
the Sun Forte Developer 7 Collection suggest a
factor of about twenty.

This is also the reason why
the trapezoidal rule, and higher rules like the Simpson Rule which follow from
a Richardson extrapolation, have been kept aside in this manuscript.

\subsection{Cubic Spline Interpolation}\label{Fm_spli.sec}
\begin{figure}
\includegraphics[bb=151 598 496 820]{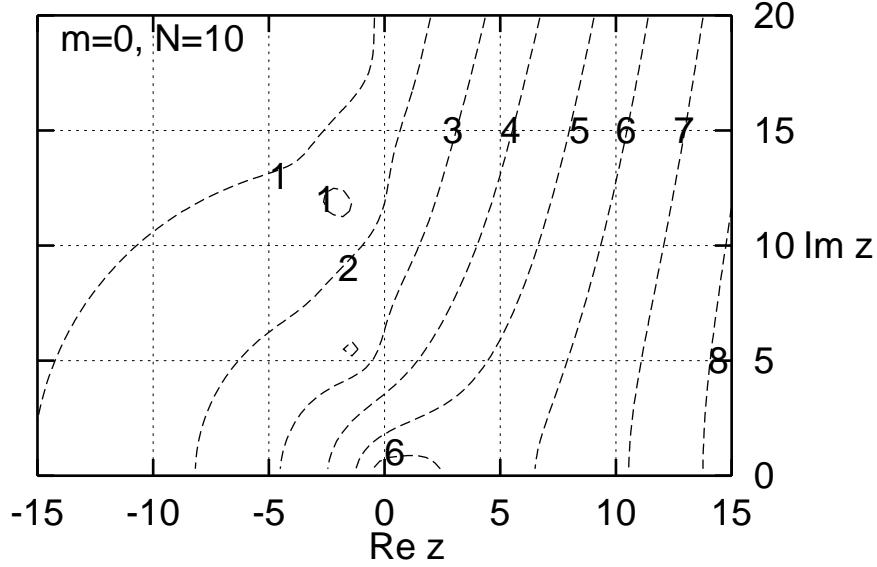}
\caption{
Contour levels of the number of valid decimal digits 
$d=1\ldots 8$ of $F_m(z)$ by the spline interpolation
(\ref{eq.Fm_spli}) for $N=10$ subintervals.
One would gain one digit throughout this $z$ region by increasing $N$ to $20$.
The accuracy for $m=3$ would be a few digits lower.
\label{fig.Fm_spli}
}
\end{figure}

The exponential of (\ref{eq.FmInt}) could be approximated by cubic splines in $N$
subintervals $[t_j,t_j+1/N]$ ($j=0,\ldots ,N-1$)
which cover the $t$-interval, to yield a sum over elementary integrals,
\begin{equation}
F_m(z)\approx\sum_{j=0}^{N-1}
\int_{t_j}^{t_{j+1}} t^{2m}
\left[ c_0+c_1 t+c_2t^2+c_3t^3\right]dt.
\label{eq.Fm_spli}
\end{equation}
In each of these intervals, the four coefficients $c_i$ are defined
by demanding that
the cubic polynomial fit and its first derivative equal the exponential
and its first derivative at both limits, $t_j$ and $t_{j+1}$:
\begin{equation}
\left(
\begin{array}{cccc}
1 & t_j & t_j^2 & t_j^3 \\
1 & t_{j+1} & t_{j+1}^2 & t_{j+1}^3 \\
0 & 1 & 2t_j & 3t_j^2 \\
0 & 1 & 2t_{j+1} & 3t_{j+1}^2 \\
\end{array}
\right)
\left(
\begin{array}{c}
c_0 \\
c_1 \\
c_2 \\
c_3 \\
\end{array}
\right)
=
\left(
\begin{array}{c}
e^{-zt_j^2} \\
e^{-zt_{j+1}^2} \\
-2zt_je^{-zt_j^2} \\
-2zt_{j+1}e^{-zt_{j+1}^2} \\
\end{array}
\right)
\label{eq.Fm_splico}
\end{equation}
The number of multiplications to solve this system of linear equations in each interval
looks prohibitive, even though one would recycle the matrix elements
and right hand sides, and even though the matrix is already close to
triangular form. One must evaluate $N$ exponentials and insert the $c_i$ into
about $2N$ polynomials of degree $2m+4$ to finalize (\ref{eq.Fm_spli}).
Actually, inserting  (\ref{eq.Fm_splico}) into (\ref{eq.Fm_spli}) yields rather
tight formulas (see \cite{Newbery} and App.\ \ref{app.quadr}), namely
\begin{equation}
F_0(z)\approx\sum_{j=0}^{N-1}
\frac{t_{j+1}-t_j}{6}
\left[
\left(3+zt_{j+1}(t_{j+1}-t_j)\right)e^{-zt_{j+1}^2}
+
\left(3-zt_j(t_{j+1}-t_j)\right)e^{-zt_j^2}
\right]
\label{eq.Fm_spli0}
\end{equation}
and
\begin{eqnarray}
F_1(z)\approx\sum_{j=0}^{N-1}
\frac{t_{j+1}-t_j}{30}
&\Big[&
\left(2zt_{j+1}^4-zt_{j+1}^2t_j^2+8t_{j+1}^2+5t_{j+1}t_j-zt_{j+1}t_j^3+2t_j^2
\right)e^{-zt_{j+1}^2} \nonumber \\
&&
 +
\left(2zt_j^4-zt_{j+1}^2t_j^2+8t_j^2+5t_{j+1}t_j-zt_{j+1}^3t_j+2t_{j+1}^2
\right)e^{-zt_j^2}
\Big] .
\label{eq.Fm_spli1}
\end{eqnarray}
The drawback formulated in Sec.\ \ref{GaussJ.sec}, however, remains:
we consider only very small $N$, and conclude from 
Fig.\ \ref{fig.Fm_spli} that this numerical expenditure is too high
to consider this method a competitive candidate.

\subsection{Salzer's Numerical Inverse Laplace Transform}
$F_m(z)$ can be written in terms of the inverse Laplace Transform of 
some function $P(m+1/2,z)$ which is loosely related to the
$\chi^2$
probability distribution \cite{Temme}:
\begin{equation}
F_m(z)=\frac{\Gamma(a)}{2z^a}P(a,z); \qquad (a\equiv m+1/2),
\label{eq.FofP}
\end{equation}
where
\begin{equation}
P(a,z)=\frac{1}{2\pi i}
\int_{c-i\infty}^{c+i\infty} e^{zs}\frac{1}{s}\frac{1}{(s+1)^a}ds; \qquad
(c>0).
\label{eq.TemmeP}
\end{equation}
The binomial expansion of $(s+1)^a$ in powers of $1/s$ and interchange
of integration and summation transforms (\ref{eq.TemmeP}) to the power
series defined by (\ref{eq.Kumm}) and (\ref{eq.KummTay}).

Substituting $zs\equiv p$ in $P(a,z)$,
(\ref{eq.TemmeP}) is in shape for Salzer's \cite{Salzer} approximation of the
kernel $1/[s(s+1)^a]$  by polynomials in $1/p$.

(\ref{eq.FofP}) and (\ref{eq.TemmeP}) assume the existence
of the Laplace Transform of $z^aF_m(z)$; this limits this proposal in general
to $\Re z>0$ as indicated by the factor
$e^z$ in (\ref{eq.KummLau}) and made explicit in \cite[p.\ 341]{MO}.
Therefore we use \cite{Temme}
\begin{equation}
P(a,z)=1-Q(a,z)=1-ae^{-z}z^a\frac{1}{2\pi i}
\int_{c-i\infty}^{c+i\infty} e^{as}\frac{1}{(as)^a}\frac{1}{z-as}ds
\label{eq.TemmeQ}
\end{equation}
to complement for $\Re z<0$. This, however, is of purely experimental
nature since (\ref{eq.TemmeQ}) depends on a contour integration
passing between the pole and the branch point, and this is not at all
accounted for in Salzer's sampling of the complex $p$-plane.

\begin{figure}
\includegraphics[bb=151 191 510 820]{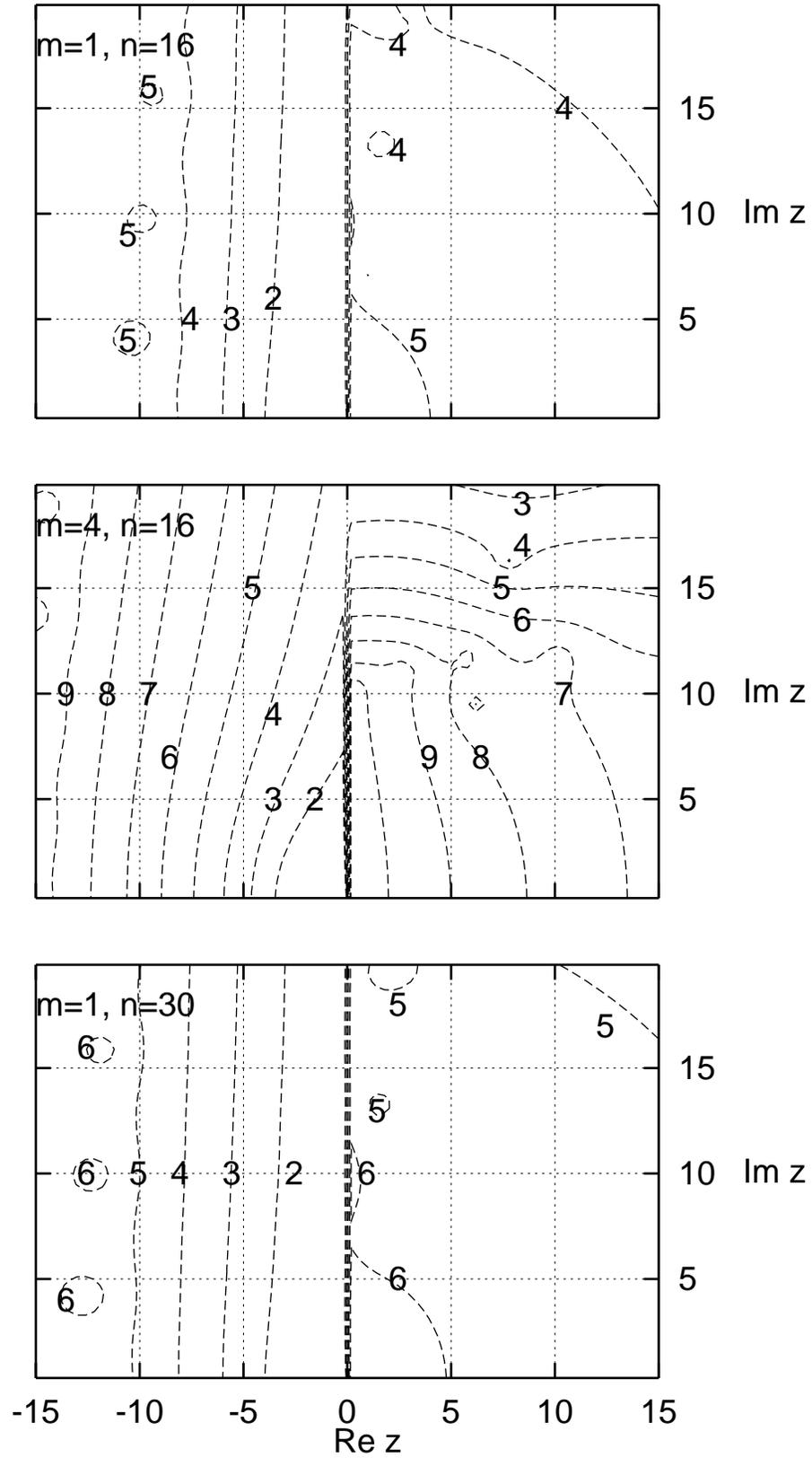}
\caption{Contour levels of the number of valid decimal digits 
$d=2\ldots 9$ of $F_m(z)$, if $n=16$ or 30 sampling points are used in
Salzer's integration of the inverse Laplace transform (\ref{eq.TemmeP}).
\label{fig.salz}
}
\end{figure}

Sample outputs of this approach are gathered in Fig.\ \ref{fig.salz}
using abscissa and weights as tabulated in Tab.\ \ref{tab.salz}.
The discontinuity in the graphs at passing the imaginary $z$-axis
is due to switching between (\ref{eq.TemmeP}) and (\ref{eq.TemmeQ}).
We see that the accuracy of $F_1(z)$ rises by about one decimal if the number $n$
of evaluations of the kernel is increased from 16 to 30, but the rise
in the case $F_4(z)$---not demonstrated in the figure---is only about a
third decimal.
The major difference to the Taylor and Lauren series methods
of sections \ref{sec.referenc}, \ref{sec.laur} and \ref{sec.kumm3} is that this
here still demands computation
of a complex valued root for each of the sampling points, and that
prediction of the accuracy both as a function of $m$ and $z$ is
complicated.

\begin{table}
\begin{tabular}{cccc}
$i$ & \multicolumn{3}{c}{$1/p_i^{(n)}$} \\
\hline
1 & 0.00837 17061 78265 71876 67547 1334&-&0.03493 88515 18794 47953 06804 3103i \\
3 & 0.01725 03391 11534 01977 09174 7117&-&0.03535 09659 30121 29557 12679 1042i \\
5 & 0.02522 57393 22044 57375 46273 8919&-&0.03348 90948 21254 18355 16245 2151i \\
7 & 0.03228 05129 34890 49074 48769 3684&-&0.02985 48245 41582 95432 08487 9898i \\
9 & 0.03823 26639 94284 05059 06029 0158&-&0.02477 02071 79697 66742 19368 1785i \\
11 & 0.04289 02498 71345 82958 76362 4358&-&0.01853 56403 81264 01685 00724 9039i \\
13 & 0.04609 15264 54310 63304 47216 7673&-&0.01146 24728 96512 75189 11368 8442i \\
15 & 0.04772 17782 62356 94180 43787 9660&-&0.00387 81037 55474 09447 26721 4718i \\
\end{tabular}

\begin{tabular}{cccc}
$i$ & \multicolumn{3}{c}{$A_i^{(n)}$} \\
\hline
1 & -(2)7.46675 12193 45759 50393 85910 48&+&(2)2.33418 71487 56825 21567 95817 62 i \\
3 & (3)2.91507 59384 65429 08402 81547 90&-&(4)6.02533 14214 97033 76429 48817 01 i \\
5 & (5)8.32343 31208 36870 55687 37470 24&+&(5)9.23999 52597 05792 07995 48624 41 i \\
7 & -(7)1.12187 25580 46183 78092 28439 34&-&(6)2.85904 20761 32552 12275 19081 14 i \\
9 & (7)5.84396 38920 01078 49663 38594 67&-&(7)1.38271 69228 73790 17106 97305 80 i \\
11 & -(8)1.53739 97073 01945 94831 81610 90&+&(8)1.17150 18184 90003 20088 57634 96 i \\
13 & (8)2.10257 24343 84449 69271 33642 99&-&(8)3.52009 23255 88077 30106 90523 00 i \\
15 & -(8)1.04572 70576 06995 39505 45148 52&+&(8)5.83415 46534 50843 20837 34946 96 i \\
\end{tabular}
\caption{Extension of Salzer's \cite{Salzer} table of 
reciprocals of the zeros of $P_n(x)$ and
Christoffel Numbers $A_i^{(n)}$ to the case $n=16$.
Numbers in parentheses are powers of $10$ as in Tab.\ \ref{tab.F0sq}.
Half of the values are not shown and follow by complex conjugation:
$1/p_i^{(n)}=1/\overline{p_{i-1}^{(n)}}$,
$A_i^{(n)}=\overline{A_{i-1}^{(n)}}$, ($i=2,4,6,\ldots,n$).
\label{tab.salz}
}
\end{table}

Finally,  we did not try to
apply this approximation to the Laplace transform of $F_m$ itself, because
sampling that kernel, put into concrete in App.\ \ref{app.lapla}\@,
demands calculation of a complex-valued
inverse trigonometric function, which is a costly
numerical task.

\clearpage

\section{Relocation to other special functions}

\subsection{Gautschi's approach to the Faddeeva Function}
\begin{figure}
\includegraphics[bb=87 398 511 826]{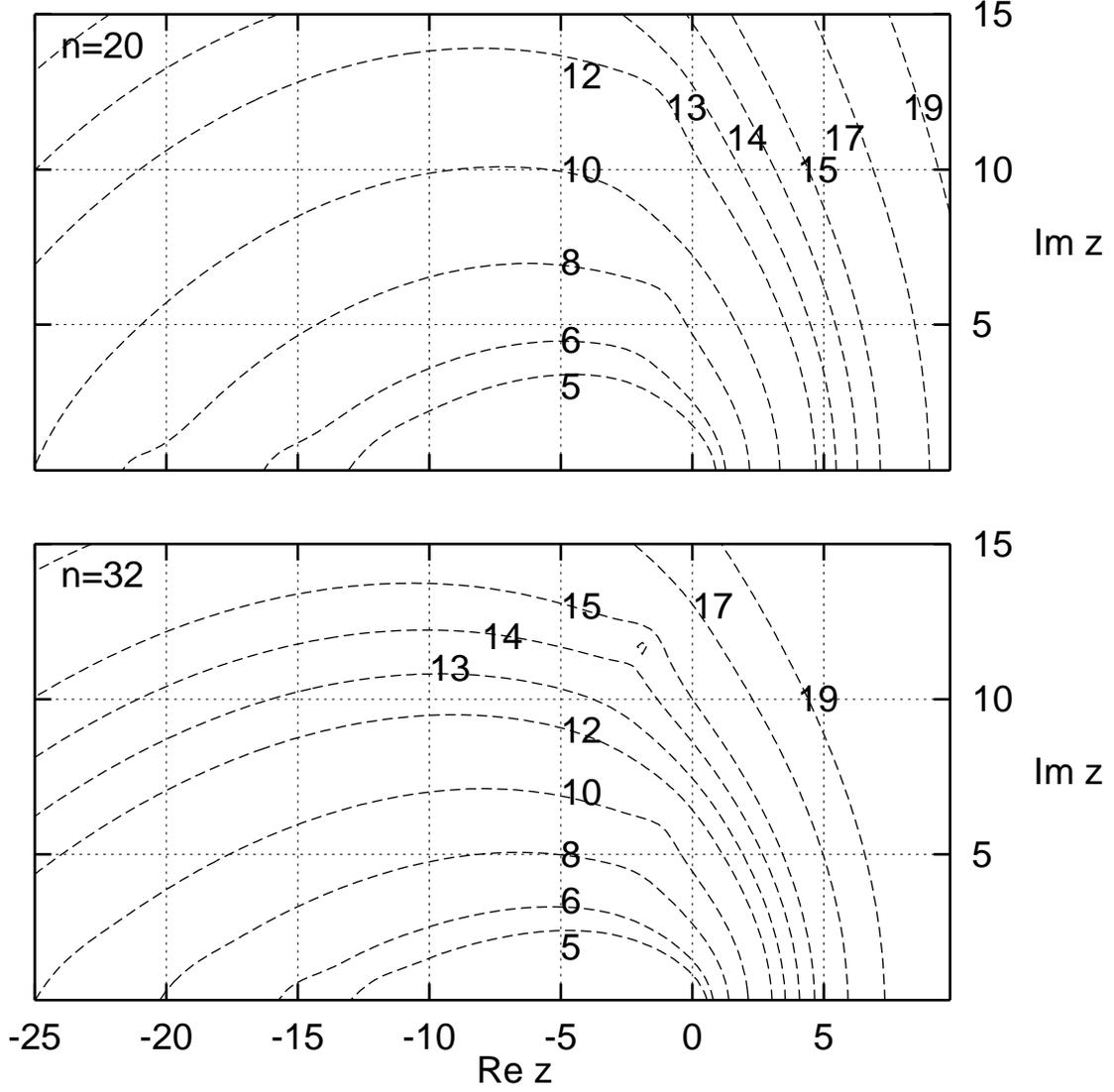}
\caption{
Contour levels of the number of valid decimal digits 
$d=5\ldots 19$ of $F_0(z)$, if the order of the Hermite Polynomial
in (\ref{eq.Fadde}) is 20 or 32.
\label{fig.gautsch}
}
\end{figure}

Following \cite{Gautschi1970,Jones} and Chapt.\ 7.1 of \cite{AS},
$F_0$ is related to the complex error function via
\begin{equation}
F_0(z)=\frac{\sqrt{\pi}}{2}\frac{\erf(\sqrt z)}{\sqrt z},
\end{equation}
and to the Faddeeva function $w$ as
\begin{equation}
F_0(z)=\frac{1}{2}\sqrt\frac{\pi}{z}\left[1-e^{-z}w(i\sqrt z)\right].
\label{eq.gaut}
\end{equation}
According to \cite[(2.14)]{Gautschi1970}\cite[(7.1.4)]{AS}\cite[(7.1.15)]{AS},
$w$ has the following representation
in terms of weights $\lambda_k^{(n)}$ and abscissae $t_k^{(n)}$ of the
$n$-point Hermite Gauss Integration (see App.\ \ref{app.Hn}):
\begin{equation}
w(z)\approx \frac{\pi}{i}\sum_{k=1}^n\frac{\lambda_k^{(n)}}{z-t_k^{(n)}}.
\label{eq.Fadde}
\end{equation}
The computation needs only about $n/2$ complex divisions
and additions, since the weights and abscissa group in symmetric pairs:
\begin{equation}
w(z)\approx \frac{2\pi}{i}z\sum_{t_k^{(n)}>0}
\frac{\lambda_k^{(n)}}{z^2-(t_k^{(n)})^2}.
\label{eq.Fadde2}
\end{equation}
Fig.\ \ref{fig.gautsch} considers the application with $n=20$ or 32.
This approach here targets the same region covered by Fig.\ \ref{fig.lauren},
either large $\Re z$ or large $\Im z$, but is obviously superior,
since it first is scalable through the choice of $n$, and at
a comparable investment into the number of complex operations
it achieves the more accurate results.

Eq.\ (2.1) in the work by Chiarella and Reichel \cite{Chiarella} looks similar
to (\ref{eq.Fadde2}): roughly speaking, the $t_k^{(n)}$ are replaced by
equidistant $nh$, and the $\lambda_k^{(n)}$ by $\exp(-n^2h^2)$.

Strand \cite{Strand} proposed a method to compute the complementary error
function of $z$, if $\Im z$ is small. It starts from the presumably known
complementary error function of $\Re z$, and therefore is
too special to be treated here.

\subsection{Expansion in Modified Spherical Bessel Functions}
\begin{figure}
\includegraphics[bb=182 236 572 701]{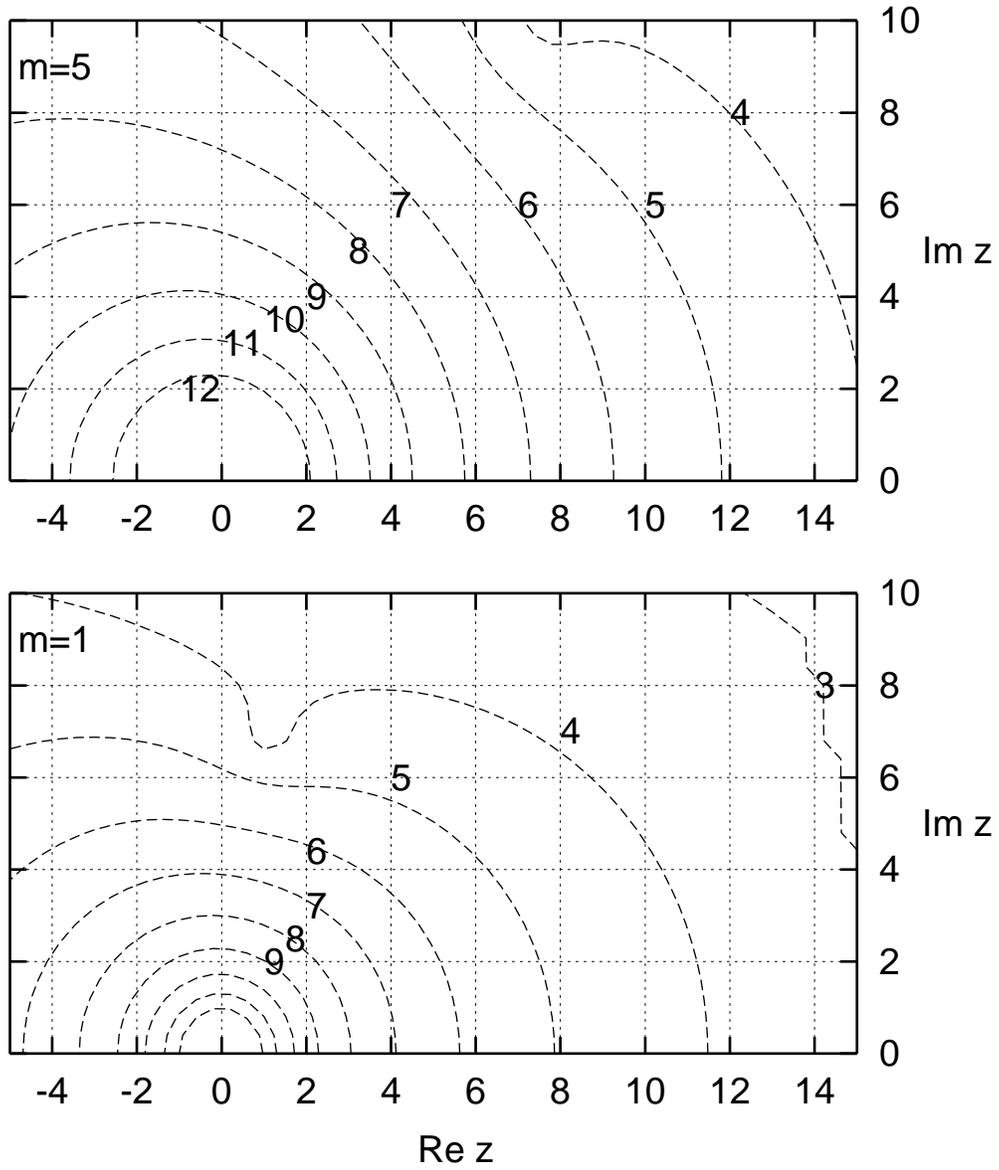}
\caption{
Contour levels of the number of valid decimals
$d=3\ldots 12$ of $F_m(z)$, if the series (\ref{eq.In})
is truncated after the term $n=7$, for $m=5$ and $m=1$.
\label{fig.In}
}
\end{figure}

The Confluent Hypergeometric Function in (\ref{eq.Kumm}) may be expanded in terms
of Modified Spherical Bessel Functions \cite[(13.3.6)]{AS}\cite{Luke1961}, which may
be rewritten with \cite[(10.2.24)]{AS} to give
\begin{eqnarray}
_1F_1(a;a+1;2z) & = &
e^z\sqrt{\frac{2\pi}{z}}\sum_{n=0}^{\infty}
(-)^n \left(\frac{1}{2}+n\right)\frac{(1-a)_n}{(1+a)_n} I_{n+1/2}(z) \\
& =&
\sum_{n=0}^{\infty}
(-)^n (2n+1)\frac{(1-a)_n}{(1+a)_n} e^z z^n 
\left( \frac {1}{z}\frac{d}{dz}\right)^n\frac{\sinh z}{z}
\label{eq.In}
\end{eqnarray}
with $a\equiv m+1/2$.
Note that a factor $(b-a-\frac{1}{2}+n)$ was missing in (13.3.6) of early editions of \cite{AS}.
The individual terms are
\begin{equation}
e^z\sqrt{\frac{2\pi}{z}}
(-)^n \left(\frac{1}{2}+n\right)\frac{(1-a)_n}{(1+a)_n} I_{n+1/2}(z)
\propto
\frac{(-)^n}{(2n-1)!!}\frac{(1-a)_n}{(1+a)_n}(z^n+z^{n+1})+O(z^{n+2}),
\label{eq.InTay}
\end{equation}
which indicates (i) that a direct implementation based on the formulas
of Tab.\ \ref{tab.In} may suffer from severe cancellation of digits
if $|z|$ is small, and
(ii) that the recurrence relations \cite[(10.2.12)]{AS} must be used in the
downward direction, for example as outlined in \cite{Mechel}.
Table \ref{tab.In}
shows that already for a small number of terms used to approximate the
series, a considerable number of complex polynomials must be computed.
Fig. \ref{fig.In} indicates that the convergence of the series is good
close to the origin of the complex plane (explained by the fact that
the lowest order terms of the Taylor series (\ref{eq.InTay})
are $O(z^n)$).
 
\begin{table}
\begin{tabular}{cc}
$n$ & $e^z z^n \left( \frac {1}{z}\frac{d}{dz}\right)^n\frac{\sinh z}{z}$ \\
\hline
   0  & $  [A-1]/[2z] $ \\
   1  & $ [z(A+1)-(A-1)]/[2z^2] $ \\
   2  & $[(z^2+3)(A-1)-3z(A+1)]/[2z^3] $ \\
   3  & $  [(z^3+15z)(A+1)-(6z^2+15)(A-1)]/[2z^4] $ \\
   4  & $  [(z^4+45z^2+105)(A-1)-(10z^3+105z)(A+1)]/[2z^5] $ \\
   5  & $  [(z^5+105z^3+945z)(A+1)-(15z^4+420z^2+945)(A-1)]/[2z^6] $ \\
   6  & $  [(z^6+210z^4+4725z^2+10395)(A-1)-(21z^5+1260z^3+10395z)(A+1)]/[2z^7] $ 
\end{tabular}
\caption{
Complexity evaluating the terms in (\ref{eq.In})
for small values of $n$, with $A\equiv \exp(2z)$.
\label{tab.In}
}
\end{table}

The assessment follows the conclusion in \cite{Schwerdtfeger1988} for a similar
expansion that was a hybrid of the Dawson and the error function for real $z$:
Since the evaluation of the Bessel Functions would be of similar complexity
as a straight-forward power series for small $|z|$, and since the
convergence is slow for large $|z|$, this type of ansatz is not competitive.

\subsection{Dijkstra's continued fractions}
\begin{figure}
\includegraphics[bb=184 200 561 824]{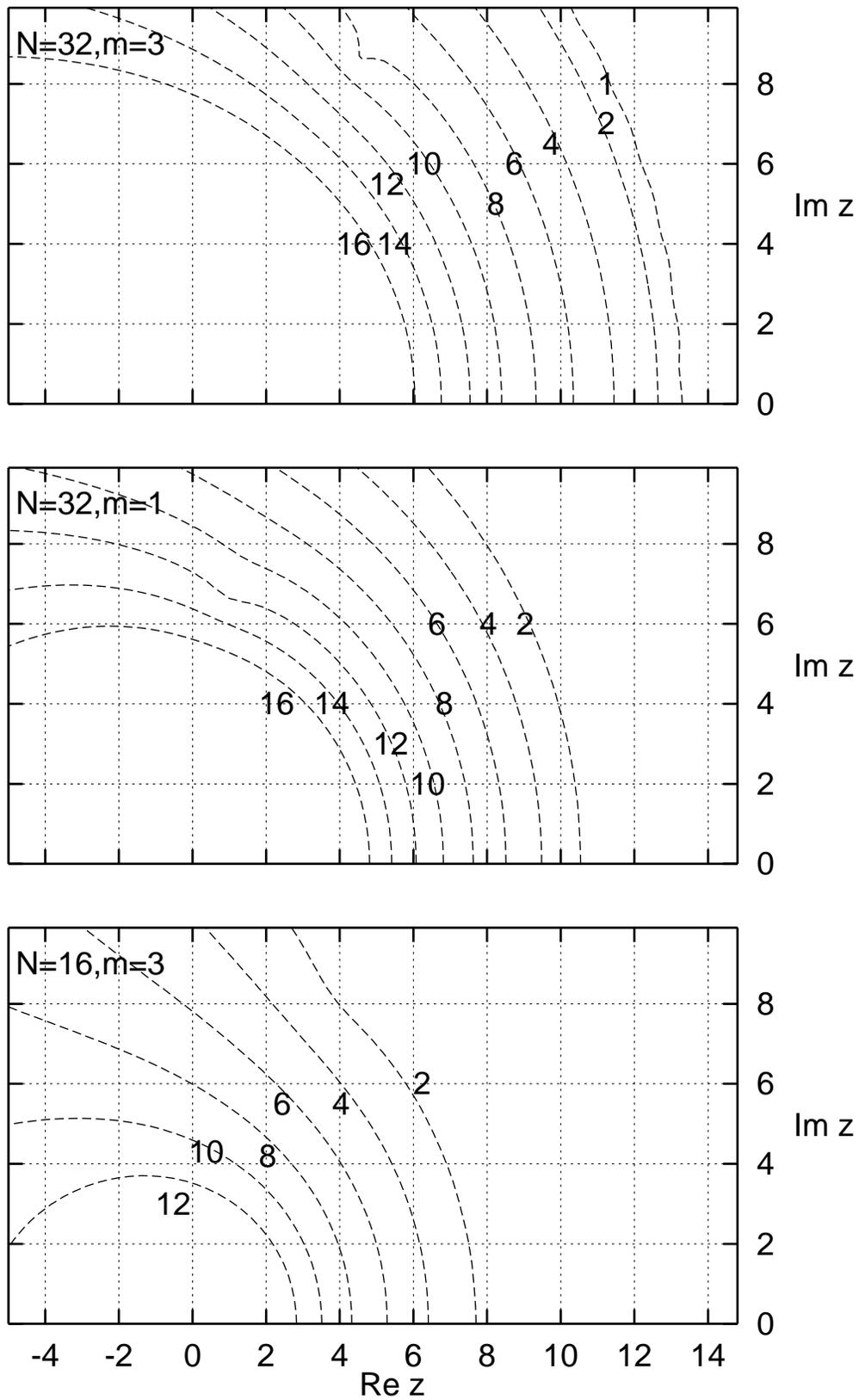}
\caption{
Contour levels of the number of valid decimal digits 
$d=1\ldots 16$ of $F_m(z)$, if Dijkstra's representation (\ref{eq.FofKaa})
is used with (\ref{eq.dijk1}) truncated at  a depths $N$ .
\label{fig.dijk}
}
\end{figure}

According to Dijkstra \cite{Dijkstra}, the constraint that (\ref{eq.KummTay})
and (\ref{eq.KummLau}) serve well only inside domains of small or large modulus
may be overcome by use of an auxiliary function $K$ that allows a suitable
continued fraction
\begin{equation}
K(a,b,z)\equiv \frac{_1F_1(a;b+1;z)}{b_1F_1(a;b;z)}
=\frac{1}{b+z-}\;\frac{z(b+1-a)}{b+1+z-}\;\cdots\frac{z(b+n-a)}{b+n+z-}\cdots
\label{eq.dijK}
\end{equation}
Application to (\ref{eq.Kumm}) with $_1F_1(a;a;-z)=e^{-z}$ proposes
\begin{equation}
F_m(z)=\frac{1}{2}e^{-z}K(a,a,-z), \qquad (a\equiv m+1/2).
\label{eq.FofKaa}
\end{equation}
In Fig.\ \ref{fig.dijk} we investigate the
readiness of the representation.
It is excellent nearby the negative real axis, staying
above $d=10$ ($m=3, N=16$), $d=13.1$ ($m=1, N=32$) and $d=16.6$
($m=3, N=32$).
Dijkstra \cite{Dijkstra} demonstrates that $K(a,b,z)$
mediates between low-$z$ and high-$z$ expansions for real positive $z$.
In so far, the sign change of $z$ in (\ref{eq.FofKaa}) means
our plot actually looks at the ``wrong'' side.
The continued fraction is terminated at the $N$th convergent.
About one decimal in accuracy has actually be gained by
deletion of the $z$ in this closing denominator like
\begin{equation}
K(a,a,z)=\frac{1}{a+z-}\;\frac{z}{a+1+z-}\cdots\frac{Nz}{a+N}.
\label{eq.dijk1}
\end{equation}

\begin{figure}
\includegraphics[bb=184 200 561 824]{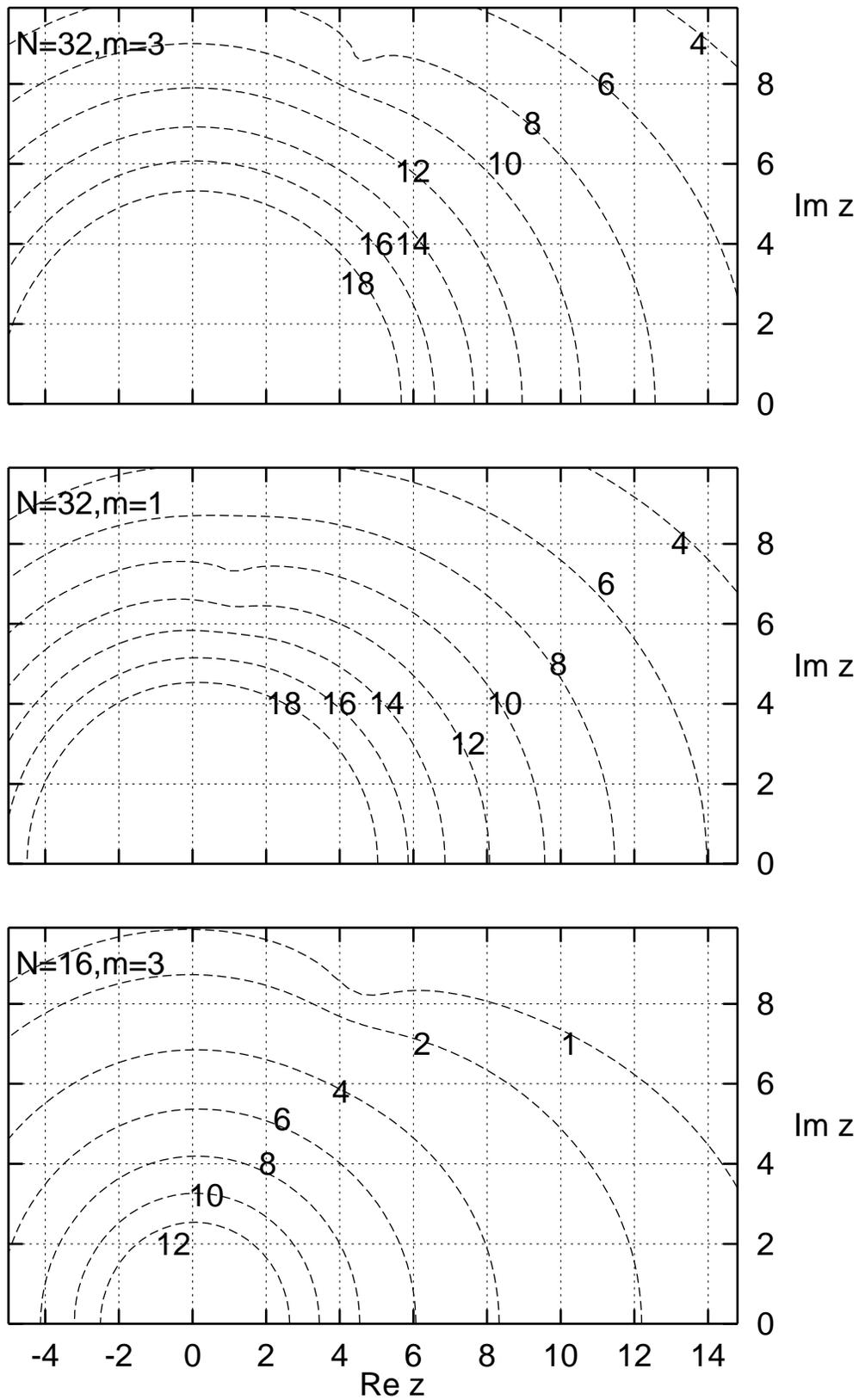}
\caption{
Contour levels of the number of valid decimal digits 
$d=1\ldots 18$ of $F_m(z)$, if Eq.\ (\ref{eq.dijk2})
is used with Dijkstra's continued fraction (\ref{eq.dijK}) truncated at a depths $N$.
\label{fig.dijk2}
}
\end{figure}

A combination of \cite[(13.4.4)]{AS}, the Kummer transformation \cite[(13.1.27)]{AS}
and the definition (\ref{eq.dijK}) enforces that both arguments of $z$
carry the same sign:
\begin{equation}
zK(1,a,z)=1-\frac{e^{-z}}{(2m-1)F_{m-1}(z)}
\label{eq.dijk2}
\end{equation}
Fig.\ \ref{fig.dijk2} shows that this indeed extends the fitness
of this Dijkstra representation into the region of $\Re z>0$,
at the expense of the fitness to compute $F_m(z)$ at $\Re z<0$.

\clearpage
\section{Summary}
Working with standard numerical methods on the fundamental integral
representation of the Incomplete Gamma Function $F_m(z)$ is
generally inefficient as it demands dense sampling (frequent
evaluation) of complex exponentials.

Continued fraction and rational function approximations are
difficult to control, because the regions of known accuracy
in the $z$-plane are complicated.

Being an analytic function of $z$, the fastest evaluation uses 
Taylor series which recall tables of the derivatives $d^nF_m(z_0)/dz_0^n$
that have been computed off-line to high precision (Sec.\ \ref{sec.kumm3}).

\clearpage
 
\appendix

\section{Example Application: Travelling Orbitals}\label{app.travo}

Let $\Psi_0({\bf r},t)=\exp(i\omega t)\varphi_0({\bf r})$ be a solution
of the Schr\"odinger equation
\begin{equation}
\left[ -\frac{\hbar^2}{2m}\nabla^2{\bf r}+V({\bf r})\right]
\Psi_0({\bf r},t)=-i\hbar \frac{\partial}{\partial t}\Psi_0({\bf r},t)
\label{eq.Schroe}
\end{equation}
associated with a stationary potential $V({\bf r})$. The 
transition to a potential that moves with constant velocity ${\bf v}$ in
the laboratory system is useful to describe electrons bound to scattering
atoms, and reads
\begin{equation}
\left[ -\frac{\hbar^2}{2m}\nabla^2{\bf r}+V({\bf r}-{\bf v}t)\right]
\Psi_{\bf v}({\bf r},t)=-i\hbar \frac{\partial}{\partial t}
\Psi_{\bf v}({\bf r},t)
\end{equation}
The solution
\begin{equation}
\Psi_{\bf v}({\bf r},t)=\exp(i\omega t)\exp(i\frac{\hbar^2k^2}{2m}t)
\exp(-i{\bf k}\cdot{\bf r})\varphi_0({\bf r}-{\bf v}t)
\end{equation}
with $\hbar{\bf k}\equiv m{\bf v}$ is generated from the solution
of (\ref{eq.Schroe}) \cite{Roos}. This Galilean transformation lets
``travelling'' orbitals $\propto \exp(-i{\bf k}\cdot({\bf r-R}))\varphi_0
({\bf r-R})$ become a natural choice for basis functions in the
laboratory coordinate system \cite{Pedersen}.
If the $\varphi_0$ are linear combinations of Gaussian Type Orbitals (GTO's),
the Coulomb Integrals may be treated with the product rule for GTO's
\cite{Boys,Zivko,Arakane,Bracken,Piccolo,Tenno,Helgaker,Obara,Dupuis1976,Lindh},
and the Gauss transform
\begin{equation}
\frac{1}{r}=\frac{2}{\sqrt\pi}\int_0^\infty e^{-r^2s^2}ds
\end{equation}
of the Coulomb potential eventually reduces the Coulomb integrals to
Incomplete Gamma Functions with complex argument $z$
\cite{Mogensen,Colle1987,Colle1988}.

\section{Laplace Representation}\label{app.lapla}

Laplace transformation of the power series (\ref{eq.Kumm}) on a term-by-term
basis yields \cite[(7.621.4)]{GR}
\begin{eqnarray}
\mathcal{L}(F_m(z))&\equiv& \int_0^{\infty}e^{-zs}F_m(z)dz
= \frac{1}{2a}\frac{1}{s}{}_2F_1(a,1;a+1;-\frac{1}{s})
= \sum_{n=1}^{\infty}\frac{(-)^{n+1}}{2(a+n-1)s^n} \\
&=&
\frac{1}{(\sqrt s)^{1-2m}}\sum_{n=1}^{\infty}\frac{(-)^{n+1}}{2m+2n-1}
\frac{1}{(\sqrt s)^{2m+2n-1}} \\
&=&
s^{a-1}\left[
(-)^{m+1}(\arctan \sqrt s-\frac{\pi}{2})+\sum_{k=1}^m
\frac{(-)^{k+1}}{(2k-1)s^{k-1/2}}
\right]
\label{eq.laplaAtan}
\end{eqnarray}
where $a\equiv m+1/2$. This result could also be derived from the Laplace
transform of the
differential equation
---which is a combination of (\ref{eq.diff}) and (\ref{eq.recurr})---:
\begin{equation}
z\frac{d}{dz}F_m(z)= -zF_{m+1}(z)=-aF_m(z)+\frac{e^{-z}}{2},
\end{equation}
\begin{equation}
-\frac{d}{ds}\left[ s (\mathcal{L}F_m)(s)-F(z=0)\right]=-a(\mathcal{L}F_m)(s)+\frac{1}{2}\frac{1}{s+1}.
\end{equation}
This inhomogeneous differential equation
\begin{equation}
s\frac{d}{ds}(\mathcal{L}F_m)(s)=(a-1)(\mathcal{L}F_m)(s)-\frac{1}{2}\frac{1}{s+1}
\end{equation}
is solved by writing down the solution of the separable homogeneous
differential equation, $\mathcal{L}F_m(s)=s^{a-1}c$, then introducing
the function $c(s)$ for the constant $c$, which leaves a simple differential
equation for $dc(s)/ds$ and
\begin{equation}
c(s)=\int \frac{1}{s^a(s+1)}ds.
\end{equation}
This is solved with the recursion \cite[(2.249)]{GR}, because $a$ is
half integer, until \cite[2.211]{GR} is applicable for the reduced remnant.

To round off this excursion: One could introduce \cite[(4.4.42)]{AS}
\begin{equation}
\arctan \alpha =\frac{\alpha}{1+\alpha^2}\left[ 1+\frac{2}{3}\frac{\alpha^2}{1+\alpha^2}
+\frac{2\cdot 4}{3\cdot 5}\left(\frac{\alpha^2}{1+\alpha^2}\right)^2
+\cdots+\frac{(2n)!!}{(2n+1)!!}\left(\frac{\alpha^2}{1+\alpha^2}\right)^n+\cdots
\right]
\end{equation}
in (\ref{eq.laplaAtan}), the simplest case of $m=0$ reading
\begin{equation}
\mathcal{L}(F_0(z))=
-\frac{1}{1+s}\left[ 1+\frac{2}{3}\frac{s}{1+s}
+\cdots+\frac{(2n)!!}{(2n+1)!!}\left(\frac{s}{1+s}\right)^n+\cdots
\right]+\frac{\pi}{2}\frac{1}{\sqrt s}.
\end{equation}
Truncation of this series after some $n$th term followed by an inverse
Laplace transform \cite[(29.2.21)]{AS}
yields $e^{-z}$ multiplied by a polynomial of degree $n$ in $z$, plus
$\sqrt{\pi/z}/2$, i.e., replaces $w(i\sqrt z)/\sqrt z$ in
(\ref{eq.gaut}) by a polynomial in $z$. Obviously, that is a poor
representation, as it forces
the polynomial in $z$ to compensate the singularity at $z=0$.

\section{Roots of Hermite Polynomials}\label{app.Hn}

These pairs of weights and abscissa of the Gauss-Hermite
quadrature can be taken from \cite[(Tab.\ 25.10)]{AS}
for some values of $n\le 20$, for $n=8$, $16$, $32$ or $64$ from \cite{Shao},
or otherwise computed with the {\tt d01bcf} routine of the NAG
library. Following on a reference by Shao et al.\ \cite{Shao} to a note
by Hofsommer \cite{Hofsommer}, the zeros of $H_n(x)$ can be refined
through a third order Newton method. Improved solutions $x^{(j+1)}$
are computed from guesses $x^{(j)}$ through \cite[(4)]{Hofsommer}, here
\begin{equation}
x^{(j+1)}=x^{(j)}-\frac{H_n(x^{(j)})}{H_n'(x^{(j)})}\left[1+x^{(j)}\frac{H_n(x^{(j)})}{H_n'(x^{(j)})}\right] .
\label{eq.Hofs}
\end{equation}
Note two sign errors in \cite[(3)]{Hofsommer}; the correct equation is
\begin{eqnarray}
\alpha&=&x-f/f'-(P+S/f')(f/f')^2 \nonumber \\
&& -\frac{1}{3}(4P^2-P'+Q+10PS/f'-S'/f'+6S^2/f'^2)(f/f')^3+O[(f/f')^4] .
\end{eqnarray}

Similar to \cite[(5.10)]{Shao} one might consider using
the terminating continued fraction
\begin{equation}
\frac{H_n'(x)}{H_n(x)}=\frac{2n}{2x-}\;\frac{2(n-1)}{2x-}\;\frac{2(n-2)}{2x-}\cdots\frac{2}{2x}
=\frac{n}{x-}\;\frac{(n-1)}{2x-}\;\frac{(n-2)}{x-}\;\frac{(n-3)}{2x-}\cdots.
\end{equation}
in (\ref{eq.Hofs}) to meet the thread of cancellation of digits.
The weights follow as \cite[(25.4.46)]{AS}
\begin{equation}
\lambda_k^{(n)} = \frac{2^{n-1}n!\sqrt{\pi}}{n^2[H_{n-1}(t_k^{(n)})]^2}.
\end{equation}

\section{``Perturbed'' Quadrature}\label{app.quadr}
Eq.\ (\ref{eq.Fm_spli0}) is an application of the Euler-Mclaurin formula
\cite[(25.4.7)]{AS}
\begin{equation}
\int_a^b f(x)dx = (b-a)\frac{f(a)+f(b)}{2}-\frac{(b-a)^2}{12}\left[ f'(b)-f'(a)\right].
\label{eq.spliGen}
\end{equation}
This is in contrast to a rule that involves the value in the midpoint
of the integration intervals \cite[(2)]{Ujevic}:
\begin{equation}
\int_a^b f(x)dx = (b-a)f\left(\frac{a+b}{2}\right)+\frac{(b-a)^2}{24}\left[ f'(b)-f'(a)\right].
\end{equation}
The evaluation of the function at the end points $a$ and $b$ in
Sect.\ \ref{Fm_spli.sec} is
cheaper than an additional evaluation in the middle of the integration interval,
because it needs two complex multiplications but no new exponentials.
So the midpoint rule and variants proposed by Hammer and Wicke
\cite{Hammer1960,Struble,Patterson,Gori} are not advantageous in our case.
A cubic spline interpolation $f(x)=\sum_{n=0}^3 c_nx^n$ induces the higher
moments
\begin{eqnarray}
\int_a^b x f(x)dx =
\frac{b-a}{60}&&\Big\{
a^2\left[2f'(b)-3f'(a)\right]
+a\left[bf'(b)+9f(b)+21f(a)+bf'(a)\right] \nonumber \\
 &&+2b^2f'(a)+21bf(b)+9bf(a)-3b^2f'(b)
\Big\},
\end{eqnarray}
\begin{eqnarray}
\int_a^b x^2 f(x)dx =
\frac{b-a}{60}&&\Big\{
a^3\left[f'(b)-2f'(a)\right]
+a^2\left[bf'(b)+4f(b)+16f(a)\right] \nonumber \\
 &&+ab\left[bf'(a)+10f(b)+10f(a)\right] \nonumber \\
&& +b^2\left[4f(a)-2bf'(b)+16f(b)+bf'(a)\right]
\Big\},
\end{eqnarray}
which establishes (\ref{eq.Fm_spli1}), and
\begin{eqnarray}
\int_a^b x^3 f(x)dx =
\frac{b-a}{420}&&\Big\{
a^4\left[4f'(b)-10f'(a)\right]
+a^3\left[5bf'(b)+90f(a)-2bf'(a)+15f(b)\right] \nonumber \\
&&+3a^2b\left[bf'(b)+13f(b)+bf'(a)+22f(a)\right] \nonumber \\
 &&+ab^2\left[66f(b)+39f(a)+5bf'(a)-2bf'(b)\right] \nonumber \\
&&+b^3\left[-10bf'(b)+4bf'(a)+15f(a)+90f(b)\right]
\Big\}.
\label{eq.x1}
\end{eqnarray}
\begin{eqnarray}
\int_a^b x^4 f(x)dx =
\frac{b-a}{840}&&\Big\{
a^5\left[5f'(b)-15f'(a)\right]
+a^4\left[7bf'(b)+150f(a)-5bf'(a)+18f(b)\right] \nonumber \\
&&+2a^3b\left[3bf'(b)+24f(b)+bf'(a)+60f(a)\right] \nonumber \\
 &&+2a^2b^2\left[42f(b)+42f(a)+3bf'(a)+bf'(b)\right] \nonumber \\
 &&+ab^3\left[-5bf'(b)+48f(a)+7bf'(a)+120f(b)\right] \nonumber \\
&& +b^4\left[-15bf'(b)+5bf'(a)+18f(a)+150f(b)\right]
\Big\}.
\label{eq.x2}
\end{eqnarray}
The generalization of (\ref{eq.spliGen}) to an integral over a quintic spline
that engages also the second derivatives (curvatures)
at the interval limits reads
\begin{equation}
\int_a^b f(x)dx = (b-a)\frac{f(a)+f(b)}{2}-\frac{(b-a)^2}{10}\left[ f'(b)-f'(a)\right]
+\frac{(b-a)^3}{120}\left[f''(b)+f''(a)\right] ,
\label{eq.x3}
\end{equation}
which simplifies to (\ref{eq.spliGen}) if $f(x)$ is any cubic polynomial.
The three formulas (\ref{eq.x1})--(\ref{eq.x3}) have not been used in this work.

\section{Gauss-Jacobi Abscissae and Weights}\label{app.gaussJ}
Formalas of weights $w_i$ and abscissae $t_i$ with (\ref{eq.gaussJ})
are given in \cite[(25.4.33)]{AS},
where $t_i$ are the zeros of Jacobi Polynomials $P_n^{(2m,0)}(1-2t)$,
and
\begin{equation}
1/w_i = \sum_{j=0}^{n-1} (2m+2j+1)[P_j^{(2m,0)}(1-2t_i)]^2 .
\end{equation}
If $m=0$, this reduces to the Gauss-Legendre quadrature,
Table 25.4 in \cite{AS}.
Table 25.8 in \cite{AS} covers the cases $m\le 2$ with $n\le 8$,
and we provide Tables \ref{tab.gaussJ2}--\ref{tab.gaussJ4} 
to cover $2m\equiv k=2$ or $4$ with $n=20$.

\begin{table}
\begin{tabular}{ccc}
$i$ & $x_i$  & $w_i$ \\
\hline
1 & (-1)0.14204 21115 93581 53319 97686 86 & (-5)0.37492 20993 33713 47725 24136 88\\
2 & (-1)0.37851 28784 95018 10290 83919 36 & (-4)0.41039 10208 73202 05549 07404 17\\
3 & (-1)0.71300 98508 12494 54654 90216 28 & (-3)0.19388 83096 17511 81076 98385 25\\ 
4 & (0)0.11385 86970 85452 85631 27170 70 & (-3)0.60704 57042 67258 25415 27237 40\\ 
5 & (0)0.16462 10853 68438 27482 84816 05 & (-2)0.14774 44441 97148 58062 83857 64\\ 
6 & (0)0.22250 74079 64513 24644 88740 54 & (-2)0.30224 89170 90737 93749 18485 98\\ 
7 & (0)0.28628 43889 84712 31925 79165 65 & (-2)0.54320 90806 12302 74500 12011 14\\ 
8 & (0)0.35459 29542 22632 43338 26791 65 & (-2)0.88135 56078 24335 53289 25273 17\\ 
9 & (0)0.42597 73418 17173 12827 24878 03 & (-1)0.13140 91925 66616 51742 10813 37\\ 
10 & (0)0.49891 61845 17151 44862 93753 46 & (-1)0.18220 50341 82686 91551 84331 53\\ 
11 & (0)0.57185 49590 66743 31538 45749 65 & (-1)0.23682 29041 42460 53980 86645 03\\ 
12 & (0)0.64323 91298 97546 89435 29449 33 & (-1)0.29002 40474 55892 85169 91635 73\\ 
13 & (0)0.71154 72878 18727 61522 80592 93 & (-1)0.33556 43145 34754 9347900271 48\\ 
14 & (0)0.77532 35806 15694 96159 06180 35 & (-1)0.36697 42047 20181 32713 25424 96\\ 
15 & (0)0.83320 87465 18990 56998 70192 99 & (-1)0.37847 39037 69223 40715 35469 55\\ 
16 & (0)0.88396 90923 44513 19505 33025 68 & (-1)0.36587 91655 32270 64986 60929 89\\ 
17 & (0)0.92652 28082 14714 71694 24334 55 & (-1)0.32734 64431 75731 48047 47964 86\\ 
18 & (0)0.95996 30995 38093 21900 36798 50 & (-1)0.26382 56442 91255 84794 51338 97\\ 
19 & (0)0.98357 79118 66012 16709 74583 86 & (-1)0.17913 75232 57385 59101 77395 77\\ 
20 & (0)0.99686 93162 59256 41043 78498 58 & (-2)0.79757 92736 27665 16852 27309 03
\end{tabular}
\caption{
Extension of \cite[Tab.\ 25.8]{AS} and \cite{Sprung} to $n=20$, $k=2$:
the abscissas and weights for the
Gaussian integration of moments,
$\int_0^1 x^kf(x)dx \approx \sum_{i=1}^n w_i f(x_i)$.
\label{tab.gaussJ2}
}
\end{table}
 
\begin{table}
\begin{tabular}{ccc}
$i$ & $x_i$  & $w_i$ \\
\hline
1 & (-1)0.28236 72218 29331 58389 89347 65 & (-7)0.17148 95670 15906 66781 55999 02 \\
2 & (-1)0.59393 81548 17751 78032 25227 62 & (-6)0.44002 56392 91105 57895 18288 23 \\
3 & (-1)0.98825 08311 50225 26221 23471 31 & (-5)0.41390 07241 19088 67017 87192 60 \\
4 & (0)0.14592 04712 27084 45276 93121 00 & (-4)0.22963 59798 81845 93019 79823 98 \\
5 & (0)0.19980 20982 26471 39701 42744 40 & (-4)0.90745 14619 04243 16233 00990 36 \\
6 & (0)0.25943 51729 24709 24862 43232 70 & (-3)0.28147 15320 39879 02977 22252 87 \\
7 & (0)0.32366 51708 61772 75707 78323 99 & (-3)0.72562 36166 25890 38003 96362 88 \\
8 & (0)0.39124 50422 79918 49074 92563 95 & (-2)0.16125 65884 61398 80639 52699 59 \\
9 & (0)0.46086 11522 78328 67423 18689 83 & (-2)0.31660 97873 98703 11170 19225 67 \\
10 & (0)0.53115 95358 51314 54342 32325 30 & (-2)0.55864 33047 22474 38067 86431 36 \\
11 & (0)0.60077 25678 99009 31814 54529 39 & (-2)0.89646 32789 96732 44746 29624 35 \\
12 & (0)0.66834 57418 12982 73134 52276 12 & (-1)0.13190 88993 08882 52772 58617 37 \\
13 & (0)0.73256 41195 17517 70122 33303 84 & (-1)0.17889 49674 56314 16340 77879 83 \\
14 & (0)0.79217 79750 00657 56854 95189 99 & (-1)0.22414 41002 15713 19313 92980 20 \\
15 & (0)0.84602 71502 09447 61169 84320 65 & (-1)0.25926 78294 51372 33371 02635 30 \\
16 & (0)0.89306 36602 60154 22796 45053 63 & (-1)0.27551 70856 12976 95932 96393 76 \\
17 & (0)0.93237 21209 85837 16946 52904 79 & (-1)0.26583 07034 97369 16858 79230 15 \\
18 & (0)0.96318 76416 89199 08774 24334 04 & (-1)0.22683 18042 40368 25551 58049 92 \\
19 & (0)0.98491 10827 62489 56330 49229 13 & (-1)0.16017 53993 12516 57464 20909 00 \\
20 & (0)0.99712 45845 24283 68493 64961 69 & (-2)0.72877 91419 97403 30297 32866 81
\end{tabular}
\caption{
Extension of \cite[Tab.\ 25.8]{AS} and Tab.\ \ref{tab.gaussJ2} to $n=20$, $k=4$.
\label{tab.gaussJ4}
}
\end{table}

\clearpage
\bibliography{math}

\end{document}